\documentclass[11pt]{article}
\usepackage[a4paper,margin=29mm]{geometry}
\usepackage[T1]{fontenc}
\usepackage[utf8]{inputenc}
\usepackage{array}
\usepackage{lmodern}
\usepackage{amsmath,amssymb,amsthm,mathtools}
\usepackage{microtype,enumitem,needspace}
\usepackage[dvipsnames]{xcolor}
\usepackage{tikz}
\usetikzlibrary{arrows.meta,calc,positioning,fit,decorations.pathreplacing}
\usepackage{aliascnt}
\usepackage[hidelinks]{hyperref}
\PassOptionsToPackage{disable}{todonotes}
\usepackage{todonotes}
\usepackage[nameinlink,noabbrev,capitalise]{cleveref}
\newif\ifshowvalidation
\newif\ifvalidationinline
\newcommand{\ValidationMarkersOff}{\showvalidationfalse}
\newcommand{\ValidationMarkersMargin}{\showvalidationtrue\validationinlinefalse}

\ValidationMarkersMargin

\newcommand{\ValidationDetails}[1]{\ifdefined\BANANAReview#1\fi}

\newcommand{\leanrepo}{https://github.com/janhubicka/lean-BANANA}
\newcommand{\leancommit}{fc0972fda0cc25ab94ab7315f6f57b7eed9127b9}
\newcommand{\successorleanrepo}{https://github.com/janhubicka/lean-successors}
\newcommand{\successorleancommit}{d9be105f73ca5a0f3f411b87d4c2419d0eb94753}

\newcommand{\ValidationLink}[6]{%
  \ifshowvalidation
    \ifvalidationinline
      \todo[inline,color=#3]{\textbf{#1.}\enspace
        \href{#5}{\textbf{Lean statement}}\enspace\nolinkurl{#6}}%
    \else
      \marginpar{\raggedright\href{#5}{%
        \textcolor{#2}{\sffamily\bfseries\scriptsize #4}}}%
    \fi
  \fi
}
\newcommand{\ValidationVerified}[2]{%
  \ValidationLink{Verified Lean statement}{green!45!black}{green!9}%
    {Lean\,\ensuremath{\checkmark}}{#1}{#2}}
\newcommand{\ValidationRelative}[2]{%
  \ValidationLink{Relative Lean verification}{green!55!black}{green!7}%
    {Lean\,rel.}{#1}{#2}}
\newcommand{\ValidationPartial}[2]{%
  \ValidationLink{Partial Lean verification}{orange!70!black}{yellow!12}%
    {Lean\,part}{#1}{#2}}

\newcommand{\leanverified}[2]{%
  \ValidationVerified{\leanrepo/blob/\leancommit/#1}{#2}}
\newcommand{\leanverifiedat}[3]{%
  \ValidationVerified{\leanrepo/blob/\leancommit/#1\#L#2}{#3}}

\newcommand{\leanverifiedsuccessorat}[3]{%
  \ValidationVerified{\successorleanrepo/blob/\successorleancommit/#1\#L#2}{#3}}

\newcommand{\leaninterfaceat}[3]{%
  \ValidationRelative{\leanrepo/blob/\leancommit/#1\#L#2}{#3}}
\newcommand{\leanpartial}[2]{%
  \ValidationPartial{\leanrepo/blob/\leancommit/#1}{#2}}
\newcommand{\leanpartialat}[3]{%
  \ValidationPartial{\leanrepo/blob/\leancommit/#1\#L#2}{#3}}

\ifdefined\BANANAPlain
  \ValidationMarkersOff
\fi
\hypersetup{
  pdftitle={Parity obstructions to precompact Ramsey expansions},
  pdfauthor={Sofia Brenner; Jan Hubicka; Colin Jahel; Jaroslav Nesetril},
  pdfsubject={Structural Ramsey theory, Fraisse classes and topological dynamics},
  pdfkeywords={structural Ramsey theory, EPPA, Fraisse classes, amenability, symplectic spaces}
}
\newcommand{\F}{\mathbb F_2}
\newcommand{\Fin}{[\mathbb Q]^{<\omega}}
\newcommand{\str}[1]{\mathbf{#1}}
\newcommand{\Age}{\operatorname{Age}}
\newcommand{\Aut}{\operatorname{Aut}}
\newcommand{\Emb}{\operatorname{Emb}}
\newcommand{\At}{\operatorname{At}}

\newcommand{\copies}[2]{\binom{#1}{#2}}

\newcommand{\BoolBANANA}{Boolean--BANANA}
\newcommand{\FolkBANANA}{Folkman--BANANA}
\newcommand{\K}{\mathcal K}
\newcommand{\Kbool}{\mathcal K_{\mathrm{bool}}}
\newcommand{\Kfolk}{\mathcal K_{\mathrm{folk}}}
\newcommand{\Ksp}{\mathcal K_{\mathrm{sp}}}

\newcommand{\Mbool}{\str M_{\mathrm{bool}}}
\newcommand{\Mboolset}{M_{\mathrm{bool}}}
\newcommand{\Mfolk}{\str M_{\mathrm{folk}}}
\newcommand{\Mfolkset}{M_{\mathrm{folk}}}
\newcommand{\Gbool}{G_{\mathrm{bool}}}
\newcommand{\Gfolk}{G_{\mathrm{folk}}}
\newcommand{\Abool}{\str A_{\mathrm{bool}}}
\newcommand{\Afolk}{\str A_{\mathrm{folk}}}
\newcommand{\Aban}[1]{\str A_{#1}}
\newcommand{\Bbool}[1]{\str B^{\mathrm{bool}}_{#1}}
\newcommand{\Bfolk}[1]{\str B^{\mathrm{folk}}_{#1}}
\newcommand{\Bban}[1]{\str B_{#1}}
\newtheorem{theorem}{Theorem}[section]
\newaliascnt{proposition}{theorem}
\newtheorem{proposition}[proposition]{Proposition}
\aliascntresetthe{proposition}
\newaliascnt{lemma}{theorem}
\newtheorem{lemma}[lemma]{Lemma}
\aliascntresetthe{lemma}
\newaliascnt{corollary}{theorem}
\newtheorem{corollary}[corollary]{Corollary}
\aliascntresetthe{corollary}
\theoremstyle{definition}
\newaliascnt{definition}{theorem}
\newtheorem{definition}[definition]{Definition}
\aliascntresetthe{definition}
\newaliascnt{example}{theorem}

\aliascntresetthe{example}
\newaliascnt{remark}{theorem}

\aliascntresetthe{remark}
\newaliascnt{question}{theorem}
\newtheorem{question}[question]{Question}
\aliascntresetthe{question}
\newaliascnt{problem}{theorem}

\aliascntresetthe{problem}
\newaliascnt{conjecture}{theorem}
\newtheorem{conjecture}[conjecture]{Conjecture}
\aliascntresetthe{conjecture}
\setlist{nosep}
\title{Parity obstructions to precompact Ramsey expansions}
\author{Sofia Brenner\thanks{Mathematical Institute, Leipzig University, Leipzig, Germany.}
\and Jan Hubi\v{c}ka\thanks{Department of Applied Mathematics, Faculty of Mathematics and Physics, Charles University, Prague, Czech Republic.}
\and Colin Jahel\thanks{Department of Mathematics, American University of Beirut, Beirut, Lebanon.}
\and Jaroslav Ne\v{s}et\v{r}il\thanks{Computer Science Institute, Faculty of Mathematics and Physics, Charles University, Prague, Czech Republic.}}

\begin{document}
\maketitle
\begin{abstract}
We give three strong amalgamation Fra\"iss\'e classes in finite languages with
function symbols whose limits are $\omega$-categorical and admit no precompact
Ramsey expansions, based on Boolean algebras with a parity map, pairs of Boolean
rings, and pairs of $\mathbb F_2$-vector spaces with a bilinear form
(\emph{BANANA}). By colouring counts whose parities are determined by the
structures modulo powers of~$2$, we obtain finite structures of infinite
Ramsey degree in each class, and likewise for every positive-dimensional
subspace in the age of the countable symplectic space over $\mathbb F_2$,
answering a question of Ivanov and Jaffrennou negatively. The BANANA class has coherent EPPA and its limit has an
amenable automorphism group with ample generics, disproving a conjecture of
Hubi\v{c}ka, Kone\v{c}n\'y, and Ne\v{s}et\v{r}il and answering a question of
Ivanov negatively. For the Boolean-algebra limit, the automorphism group admits
a metrisable minimal flow with meagre orbits and continuum many nonconjugate
closed maximal amenable subgroups.
\end{abstract}

\tableofcontents

\section{Introduction}\label{sec:intro}

A class $\mathcal C$ of finite structures is a \emph{Ramsey class} if for
every $\str A,\str B\in\mathcal C$ and $r\ge1$ there is $\str C\in\mathcal C$
such that, for every $r$-colouring of the embeddings of $\str A$ into $\str C$,
there exists an embedding $f\colon\str B\to\str C$ for which all composite
embeddings $f\circ e$, with $e\colon\str A\to\str B$ an embedding, have the
same colour. More generally, the \emph{Ramsey degree} of $\str A$ in
$\mathcal C$ is the least $t\ge1$ such that for every $\str B\in\mathcal C$
and $r\ge1$ there is $\str C\in\mathcal C$ such that, for every $r$-colouring
of the embeddings of $\str A$ into $\str C$, there exists an embedding
$f\colon\str B\to\str C$ for which the composites $f\circ e$ use at most $t$
colours, or infinity if no such $t$ exists (in the literature, Ramsey classes
and degrees are sometimes defined with respect to copies---substructures
isomorphic to $\str A$---rather than embeddings; the two notions are not
equivalent when $\str A$ has nontrivial automorphisms, and we call the copy
variant the \emph{copy Ramsey degree}). An \emph{expansion} enriches the structures of $\mathcal C$ by
additional relations, functions or constants. A hereditary,
reasonable expansion class with the Ramsey property is a \emph{Ramsey
expansion} (see \Cref{sec:part1-preliminaries}). It is \emph{precompact} if
copies of each $\str A\in\mathcal C$ generate only finitely many expanded
substructures up to isomorphism (for expansions by relations, this simply
means that each $\str A$ admits only finitely many expansions up to
isomorphism).

Structural Ramsey theory emerged in the 1970s. Leeb isolated the notion of
an $\str A$-Ramsey property for categories, developed in the work of
Graham--Leeb--Rothschild \cite{GLR,LeebPascal73}. In 1977,
Ne\v set\v ril--R\"odl proved the Ramsey theorem for finite ordered
relational structures, even with prescribed irreducible (clique-like)
substructures forbidden \cite{NRBull77,NR}; independently,
Abramson--Harrington proved the case without forbidden substructures
in 1978 \cite{AH78}. Further Ramsey classes and
methods were developed throughout the 1980s (see the recent survey~\cite{HK} for references). The 2005
Kechris--Pestov--Todorcevic correspondence connects these finite-colouring
problems, via Ramsey expansions, with the dynamics of automorphism groups
\cite{KPT}. The correspondence revived interest in Ramsey classes
and led to a systematic framework based on precompact Ramsey expansions
with the expansion property, developed further by Nguyen Van Th\'e and Zucker
\cite{NVTPrecompact13,ZuckerDynamics16}; see also Zucker's work on the generic
point problem \cite{ZuckerGenericPoint18}. This renewed activity also led to
general methods for constructing Ramsey classes, notably the
Hubi\v cka--Ne\v set\v ril framework for classes with closures (functions) and forbidden
homomorphisms \cite{HNAllThoseRamsey}. Precise definitions,
including the hereditary Fra\"iss\'e-age convention used here, are given in
\Cref{sec:part1-preliminaries}.%

In 1989, and in a broader formulation in 2005, Ne\v set\v ril
observed that a hereditary Ramsey class with the joint embedding property has
the amalgamation property \cite{Nesetril89,Nesetril05,HK}. Another
strengthening of amalgamation is the \emph{extension property for partial
automorphisms} (\emph{EPPA}), also called the \emph{Hrushovski property}.
A class $\mathcal C$ has EPPA if for every $\str A\in\mathcal C$ there is
$\str B\in\mathcal C$ containing $\str A$ such that every isomorphism
between two substructures of $\str A$ extends to an automorphism of
$\str B$ \cite{HKN22}.

EPPA developed independently of structural Ramsey theory. Hrushovski proved
it for finite graphs in 1992 \cite{Hrushovski92}. 
In 1993, it was a central tool in the proof of the small index conjecture by
Hodges--Hodkinson--Lascar--Shelah \cite{HHLS93}.
Herwig gave a generalisation to a broad class of finite structures in 1995 \cite{Herwig95}. 
This was followed by the
seminal Herwig--Lascar theorem in 2000 \cite{HerwigLascar00}. In 2007,
Kechris--Rosendal introduced the framework of ample generics and connected it
to amalgamation properties of finite systems of partial automorphisms
\cite[Section~6]{KR}. Macpherson's 2011 survey treated EPPA and structural
Ramsey theory in the same homogeneous-structure framework and explicitly
asked for a direct connection between them
\cite[Problem~6.5.4]{MacphersonSurvey}; see also Siniora's 2017
thesis \cite{SinioraThesis}. In 2020, Siniora--Solecki connected coherent
EPPA with dynamical properties of automorphism groups, obtaining dense locally finite
subgroups (which imply amenability) and ample generics for broad
free-amalgamation examples \cite{SS,PestovDynamics}. The two theories
subsequently began to converge. Evans, Hubi\v cka and Ne\v set\v ril
developed common structural methods for both Ramsey and EPPA in 2021 \cite{EHN21}, and
the 2022 paper~\cite{HKN22} by Hubi\v cka, Kone\v cn\'y, and Ne\v set\v ril gave a systematic strengthening
of the Herwig--Lascar theory. 

An important question in the area, also formulated as a conjecture by
Nguyen Van Th\'e \cite[Conjecture~1]{NVTsurvey}, was the following:
\begin{question}[Bodirsky--Pinsker--Tsankov~\cite{BPT}; see also~\cite{MNT}]\label{q:ramsey_precompact}
	Does the age of every $\omega$-categorical structure have a precompact Ramsey expansion?
\end{question}
Evans answered it negatively by using a
Hrushovski predimension construction to obtain the first countable
$\omega$-categorical structure with no precompact Ramsey expansion (also giving
the first Roelcke-precompact Polish group with nonmetrisable universal minimal
flow, answering~\cite[Question~1.5]{MNT}). The
argument was circulated in his private note and published with Hubi\v cka
and Ne\v set\v ril in 2019 \cite{EvansHrushovskiPrivate,EHN}. The
corresponding age, however, does not have EPPA \cite{EHN}. %

In this paper, we construct three classes that exhibit increasingly flexible embeddings, have $\omega$-categorical Fra\"iss\'e limits and admit no precompact Ramsey expansions.
The class \emph{\BoolBANANA} consists of finite Boolean algebras with an additive parity
map that takes the value zero on the greatest element. \emph{\FolkBANANA} drops the
requirement that embeddings preserve the greatest element; to prepare the
comparison with BANANA, it is formulated with two Boolean rings and a map to
$\mathbb F_2$ additive in each variable. \emph{BANANA} retains only the additive
groups on the two sides and their bilinear pairing. 
(The name BANANA abbreviates \emph{Bipartite Additive
Neighbourhood Age, Naturally Amalgamated}: neighbourhoods add by symmetric
difference, and new cross-pairings in an amalgam may be chosen freely.)
To our knowledge, it is the first class with EPPA and an $\omega$-categorical Fra\"iss\'e limit that
has no precompact Ramsey expansion. 

\Needspace{13\baselineskip}
\begin{theorem}\label{thm:overview}
The three classes \BoolBANANA, \FolkBANANA, and BANANA are strong amalgamation Fra\"iss\'e classes in finite
languages with function symbols, and their limits are $\omega$-categorical.
None has a precompact Ramsey expansion.
For BANANA, writing a member as $\str A=(L_A,R_A,\beta_{\str A})$, we have in addition:
\begin{enumerate}[label=\textup{(\roman*)}]
\item the class has coherent EPPA;
\item the automorphism group of its limit is amenable and has ample generics;
\item the Ramsey degree of $\str A$ is $|\Aut(\str A)|$ (its copy Ramsey
degree is $1$) if $L_A=0$ or $R_A=0$, and infinite otherwise;
\item the automorphism group of the limit has no co-precompact extremely
amenable subgroup and has a nonmetrisable universal minimal flow. Equivalently,
the limit has no $\omega$-categorical expansion in a countable language on the
same underlying set with extremely amenable automorphism group.
\end{enumerate}
\BoolBANANA{} and \FolkBANANA{} do not have EPPA, and the automorphism
groups of their limits are not amenable.
\end{theorem}

All three structures are defined precisely in \Cref{sec:class}. Their embeddings are described, respectively,
by partitions, disjoint unions and injective linear maps.

The structures considered are homogeneous in finite languages with function symbols, not finite relational languages. 
For BANANA, the limit has the particularly simple concrete presentation
\[
 L\simeq R\simeq[\mathbb Q]^{<\omega},\qquad
 E(X,Y)\ \Longleftrightarrow\ |X\cap Y|\text{ is odd},
\]
where $[\mathbb Q]^{<\omega}$ denotes the set of finite subsets of
$\mathbb Q$, and $L$ and $R$ are disjoint copies of it. On each copy,
addition is given by symmetric difference.

As a model-theoretic side consequence (\Cref{sec:definable-ramsey}), the
first-order theory of the BANANA limit (unlike its automorphism group) is
extremely amenable in the sense of
Hrushovski--Krupi\'nski--Pillay~\cite{HKPFirstOrder}, equivalently it has the
definable elementary embedding Ramsey property of
Krupi\'nski--Lee--Moconja~\cite{KLMRamseyDynamics}. Thus definable finite
colourings can always be made monochromatic, whereas for arbitrary finite
colourings every BANANA structure $\str A$ with $L_A\ne0$ and $R_A\ne0$ has
infinite Ramsey degree.

\paragraph{Questions answered.}\label{sec:questions-answered}
Hubi\v cka--Kone\v cn\'y--Ne\v set\v ril formulated the following
conjecture \cite[Conjecture~1.7]{HKN22}, repeated in
\cite[Conjecture~8.1]{HK}. We use its Fra\"iss\'e-age/KPT formulation,
so \emph{Ramsey expansion} has the hereditary meaning fixed in
\Cref{sec:part1-preliminaries}.
\begin{conjecture}\label{conj:eppa-ramsey}
Every strong amalgamation class with EPPA has a precompact Ramsey expansion.
\end{conjecture}
The BANANA class is a counterexample to this conjecture.%

Ivanov asked its amenable, $\omega$-categorical counterpart \cite{Ivanov}.
Evans--Hubi\v cka--Ne\v set\v ril later highlighted this question in
the structural-Ramsey setting \cite[Question~7.3]{EHN}.
\begin{question}[Ivanov]\label{q:ivanov}\label{q:ivanov-group}
Does every countable $\omega$-categorical structure with amenable
automorphism group have an $\omega$-categorical expansion in a countable
language whose automorphism group is extremely amenable? Equivalently, must a closed
amenable oligomorphic permutation group on a countable set contain a
co-precompact extremely amenable subgroup?
\end{question}
The BANANA class is a counterexample to this question as well, even with coherent EPPA.
The proof of \Cref{thm:overview} occupies
\Cref{sec:class,sec:eppa,sec:no-precompact}.

The method presented here extends beyond BANANA structures and can be generalised
to bilinear forms. In particular, we can negatively answer the following question
\cite[Question~A]{IJ}.
\begin{question}[Ivanov--Jaffrennou]\label{q:IJ-main}
Do the ages of countably infinite-dimensional symplectic, orthogonal or
unitary spaces over finite fields have finite Ramsey degrees?
\end{question}
\begin{theorem}\label{thm:symplectic-overview}
In the age of the countably infinite-dimensional symplectic space over
$\mathbb F_2$, every positive-dimensional member has infinite Ramsey
degree. The zero-dimensional member has degree $1$.
\end{theorem}
\leanverifiedat{BANANA/NonPrecompact/RamseyDegreeClassifications.lean}{28}{symplectic embedding Ramsey degree classification}
Here the age includes degenerate subspaces. Already a one-dimensional
subspace has infinite Ramsey degree. 
Thus \Cref{thm:symplectic-overview} strengthens the failure of the Ramsey property proved in \cite{IJ}.
The proof is given in \Cref{sec:symplectic-section}.

\paragraph{Further \BoolBANANA{} dynamics.}
The class \BoolBANANA{} also gives two complementary dynamical conclusions, proved
in \Cref{sec:prebanana-dynamics}. Its automorphism group has an explicit compact
metrisable minimal flow all of whose orbits are meagre (giving an elementary
counterpart to~\cite[Theorem~5.2]{EHN}), and it has continuum many pairwise
nonconjugate closed maximal amenable subgroups.

\paragraph{Comparison with earlier examples.}\label{sec:counterexample-landscape}
\begin{table}[htbp]
\centering\scriptsize
\setlength{\tabcolsep}{2pt}
\renewcommand{\arraystretch}{1.16}
\begin{tabular}{@{}lccccc@{}}
\hline
Example & \shortstack{Strong\\amalgamation} & $\omega$-cat.
 & \shortstack{Homogeneous\\in finite\\language}
 & \shortstack{Amenable\\group} & EPPA\\
\hline
$(\mathbb Z,d)$ & No & No & Yes & Yes & No\\
Tree-like metric ages~\cite{MetricallyHomogeneous} & No & No & Yes & --- & No\\
Wa\.zewski $W_P$ ($P$ infinite)~\cite{KwiatkowskaDendrites} & n/a & No & No & --- & No\\
EHN sparse (non-$\omega$-cat.)~\cite{EHN} & Yes$^*$ & No & No & No & No\\
EHN sparse ($\omega$-cat.)~\cite{EHN} & Yes$^*$ & Yes & No & No & No\\
\BoolBANANA & Yes & Yes & Yes & No & No\\
\FolkBANANA & Yes & Yes & Yes & No & No\\
BANANA & Yes & Yes & Yes & Yes & Coherent\\
\hline
\end{tabular}
\caption{Explicit structural examples with no precompact Ramsey expansion. The asterisk in the EHN rows refers to strong embeddings.}
\label{tab:examples}
\end{table}

Table~\ref{tab:examples} places the three classes alongside several earlier
examples of classes without precompact Ramsey expansions.
The fourth column records whether the limit can be presented as a homogeneous
structure in a finite language of relations and functions without changing its
automorphism group.

The metric space $(\mathbb Z,d)$ of the integers already gives an elementary
example without a precompact Ramsey expansion, %
but its isometry group is discrete and pairs have
infinitely many orbits. The tree-like metrically homogeneous graphs of
Aranda--Bradley-Williams--Hubi\v cka--Karamanlis--Kompatscher--Kone\v cn\'y--Pawliuk
generalise this metric phenomenon: the ages of their associated metric
spaces have no precompact Ramsey expansion and no EPPA
\cite[Theorems~1.1 and~1.2]{MetricallyHomogeneous}.

Kwiatkowska's generalised Wa\.zewski dendrites give a different and
particularly relevant comparison \cite{KwiatkowskaDendrites}. Her finite
tree family $\mathcal F_P$ has joint embedding and amalgamation but is not
hereditary under arbitrary substructures. 
She works with the corresponding
Fra\"iss\'e--HP variant, using a restricted hereditary notion. When $P$ is
infinite, she proves that a fixed two-vertex member has infinite Ramsey
degree in this embedding-based framework
\cite[Theorem~23]{KwiatkowskaDendrites}. 
She also constructs a reasonable
expansion $\mathcal F_P^*$ with the corresponding restricted hereditary,
expansion and Ramsey properties. For infinite $P$, the infinite-degree and
universal-minimal-flow results force this expansion to be nonprecompact
\cite[Theorems~13 and~23, Corollary~24]{KwiatkowskaDendrites}.
Equivalently, the homeomorphism group of $W_P$ has a nonmetrisable universal
minimal flow, in fact with a comeagre orbit.

The sparse-graph work of Evans--Hubi\v cka--Ne\v set\v ril gives the $\omega$-categorical example %
discussed earlier, whose automorphism group is not
amenable \cite{EHN}. There, one first passes to closure reducts encoded by
partial functions and then obtains precompact Ramsey expansions of those
enriched classes
\cite[Definition~6.4, Lemma~6.6 and Theorem~6.7]{EHN}. BANANA differs from
all these comparison examples in that it is already $\omega$-categorical and
amenable, has coherent EPPA and strong amalgamation, yet its Ramsey expansion
is nonprecompact.

More recent work of Basso--Tsankov~\cite{BassoTsankovKaleidoscopic} extends the dendrite phenomenon at the
level of topological groups. 
Their kaleidoscopic construction gives a large
class of non-Archimedean Polish groups---that is, groups with a neighbourhood
basis at the identity consisting of open subgroups---with nonmetrisable
universal minimal flows, often still with a comeagre orbit.
There are also general dynamical ways of producing
nonprecompact examples from groups with nonmetrisable universal minimal
flows. 
We do not list these in Table~\ref{tab:examples}, since the table is intended to record the explicit 
examples closest to the present paper.

\paragraph{Open questions.}
Several related questions remain open. The most notable is the following.
\begin{question}[{Bodirsky--Pinsker--Tsankov~\cite[Problem~27]{BPT}}]\label{q:ramsey_finite}
	Does every amalgamation class in a finite relational language have a Ramsey expansion in a finite relational language?
\end{question}
While BANANA can also be based on a class of $C$-relations given by the Milliken tree theorem (a Milliken--BANANA),
the parity information imposes enough structure that the resulting structure is not homogenisable in a finite relational language.

Finally, Kwiatkowska--Malicki asked whether a non-Archimedean Polish
group can be simultaneously extremely amenable and have ample generics
\cite{KwiatkowskaMalicki}. BANANA has ample generics but is not
itself extremely amenable, so it does not answer this question.

\paragraph{How the obstruction was found.}
The route to the present examples began with a different project. At the
BIRS workshop \emph{Infinite Structural Ramsey Theory}, held in Banff on
23--28 November 2025, Colin Jahel and Sofia Brenner began investigating Ramsey expansions arising from amalgamation
classes of finite groups. Together with Jan Hubi\v cka, they first analysed finite
abelian groups and then tried to pass to groups of nilpotency class at most
two. 
This project is developed in the companion manuscript~\cite{BrennerHubickaJahelWIP}.

A natural class of groups to start with was Baudisch's Fra\"iss\'e class of finite groups of nilpotency class $2$ and odd prime exponent,
with a named central subgroup containing the commutator subgroup~\cite{Baudisch01}.
Via the Lazard correspondence and a splitting, its Ramsey problem leads to
the analogous questions for alternating bilinear systems
\cite{BrennerHubickaJahelWIP}.
Closely related polar and symplectic geometries also play a central role in the
structure theory of Cherlin--Hrushovski \cite{CherlinHrushovski03}.  
The bilinear Ramsey question itself had been raised earlier: following discussions
in Edinburgh, David Evans asked Jan Hubi\v cka in 2021 whether canonically ordered
vector spaces with bilinear forms might form a Ramsey class
\cite{EvansBilinear2021}.

Our first approach was therefore to prove a Ramsey theorem
for bilinear systems equipped with admissible orders.
The 2025 result of
Ivanov--Jaffrennou that the unexpanded symplectic classes are not Ramsey did
not rule out this possibility \cite{IJ}.  When the direct attempts did not
produce a proof, we followed the systematic programme proposed in
\cite{HubickaOberwolfach25}: rather than guessing an expansion, inspect the
tree of finite types of an enumeration and try to read the missing Ramsey
information from it.  At this stage, for the first time in our work, we also
experimented with AI-assisted tools to enumerate small type configurations,
compare candidate records and test proposed invariants.  These
experiments were exploratory. 
The mathematical statements used here were
subsequently proved and checked independently.  The same analysis produced,
in odd characteristic, a fixed bilinear system of infinite Ramsey
degree and hence an obstruction to precompact Ramsey expansions in the same
Fra\"iss\'e-age sense. 
The obstruction also lifts back to the Baudisch class
\cite{BrennerHubickaJahelWIP}.

A tree-of-types analysis exposes a finite structure one coordinate or
refinement at a time, recording at each node the type visible so far. A
typical Ramsey argument passes to a substructure or subtree on which extension
steps carrying no essential new type information can be omitted. In our
examples, however, an apparently redundant step may
still change a parity value. 
For instance, one more common coordinate
switches the parity of $|X\cap Y|$, and such events cannot be discarded. 
Surprisingly, this yields persistent colourings that colour the
sets by their sizes, showing that this is not just a limitation of the method,
but an actual obstacle to the class being Ramsey.

Clearly, the existence of many possible integer values is not enough. 
However, the essential step to obtain the lower bound is to prove that many residues occur inside \emph{every} copy of a suitable finite target. 

\paragraph{Notation for the main classes.}
The notation used for the classes studied in the paper is summarised below.
\Needspace{8\baselineskip}
\begin{center}
\small
\renewcommand{\arraystretch}{1.12}
\begin{tabular}{@{}lll@{}}
\hline
Class & Notation & Limit or ambient structure\\
\hline
\BoolBANANA & $\Kbool$ & $\Mbool$\\
\FolkBANANA & $\Kfolk$ & $\Mfolk$\\
BANANA & $\K$ & $\str M$\\
symplectic age over $\F$ & $\Ksp$ & countable symplectic space\\
positive $S$-mass \BoolBANANA & $\Kbool^{\rho,S}$ & $(\Mbool,\rho_S)$\\
\hline
\end{tabular}
\end{center}

This paper is structured as follows: in \Cref{sec:part1-preliminaries}, we fix the terminology. Sections~\ref{sec:class} and~\ref{sec:eppa} give
the Fra\"iss\'e, EPPA and amenability constructions. 
In Section~\ref{sec:no-precompact}, we describe the finite Ramsey obstructions and their dynamical consequence. 
In Section~\ref{sec:symplectic}, we give the symplectic application. 
Finally, Section~\ref{sec:prebanana-dynamics} records two further dynamical
features of \BoolBANANA.

\section{Preliminaries}\label{sec:part1-preliminaries}

We follow the terminology of \cite{HK}. 
Unless stated otherwise, all structures are finite. 
Structures are denoted by boldface letters ($\str A,\str B,\str M,\ldots$) and
their underlying sets by the corresponding lightface letters
($A,B,M,\ldots$). We freely use finite many-sorted languages.
An \emph{embedding} is a
sort-preserving injection which preserves functions and constants and preserves
and reflects relations. Substructures are closed under the functions and contain
the named constants; the substructure \emph{generated} by a subset is the
smallest substructure containing it.

A class of structures $\mathcal C$ has the \emph{hereditary property} if it is closed under taking finite substructures. 
It has the \emph{joint embedding property} if for any two members of $\mathcal C$, there is a member of $\mathcal C$ that embeds both. 
It has the \emph{amalgamation property} if any two embeddings with common domain can be completed
to a commuting square. 
An amalgam is \emph{strong} if the two images intersect exactly in the image of the common domain. A \emph{Fra\"iss\'e class} is a
hereditary class with joint embedding and amalgamation, countably many
isomorphism types and arbitrarily large members. 
Its \emph{Fra\"iss\'e limit} (unique up to isomorphism) is the countable locally finite homogeneous structure with that age. Here \emph{locally finite}
means that finitely generated substructures are finite, and \emph{homogeneous}
means that every isomorphism between finite substructures extends to an
automorphism. 
The \emph{age} $\Age(\str M)$ of $\str M$ is the class of finite structures embeddable in~$\str M$.

A permutation group is \emph{oligomorphic} if it has finitely many orbits on
$n$-tuples for every~$n$. 
By the Ryll--Nardzewski theorem \cite[Theorem~3.1.2]{MacphersonSurvey},
for a countably infinite structure in a countable language, oligomorphicity
of its automorphism group is equivalent to $\omega$-categoricity, i.e., that its first-order theory has a unique countable model up to isomorphism.

\paragraph{Ramsey degrees and expansions.}
Write $\Emb(\str A,\str B)$ for the set of embeddings of $\str A$ into $\str B$ and $\copies{\str B}{\str A}$ for the set of their images, also called \emph{copies} of $\str A$ in $\str B$.
For $r,t\ge1$, write $\str C\longrightarrow(\str B)^{\str A}_{r,t}$ if every $r$-colouring of $\Emb(\str A,\str C)$ takes at most $t$ colours on $\{f\circ e:e\in\Emb(\str A,\str B)\}$ for some $f\in\Emb(\str B,\str C)$, and abbreviate $\str C\longrightarrow(\str B)^{\str A}_{r,1}$ to $\str C\longrightarrow(\str B)^{\str A}_r$.
A class $\mathcal C$ of finite structures is a \emph{Ramsey class} (or has the \emph{Ramsey property}, also called the \emph{embedding Ramsey property}) if for every $\str A,\str B\in\mathcal C$ and every $r\ge1$ there is $\str C\in\mathcal C$ with $\str C\longrightarrow(\str B)^{\str A}_r$.
The \emph{Ramsey degree} $t_{\mathcal C}(\str A)$ of $\str A\in\mathcal C$ is the least $t\ge1$ such that for every $\str B\in\mathcal C$ and $r\ge1$ there is $\str C\in\mathcal C$ with $\str C\longrightarrow(\str B)^{\str A}_{r,t}$, or $\infty$ if no such $t$ exists.
Similarly, for colourings of copies rather than embeddings, write $\str C\longrightarrow_{\rm cp}(\str B)^{\str A}_{r,t}$ if every
$r$-colouring of $\copies{\str C}{\str A}$ uses at most $t$ colours on $\copies{\str B'}{\str A}$ for some $\str B'\in\copies{\str C}{\str B}$. The \emph{copy Ramsey degree} $t^{\rm cp}_{\mathcal C}(\str A)$ is
the least $t\ge1$ such that for every $\str B\in\mathcal C$ and $r\ge1$ there is
$\str C\in\mathcal C$ with $\str C\longrightarrow_{\rm cp}(\str B)^{\str A}_{r,t}$, or $\infty$ if no such $t$ exists.
For every $\str A\in\mathcal C$ we have $t_{\mathcal C}(\str A)=|\Aut(\str A)|\,t^{\rm cp}_{\mathcal C}(\str A)$ \cite{HK}: fixing an enumeration of $A$ and a linear order on $C$ orders lexicographically the $m=|\Aut(\str A)|$ embeddings of $\str A$ onto each $\str A'\in\copies{\str C}{\str A}$, so an $r$-colouring of $\Emb(\str A,\str C)$ induces an $r^m$-colouring of $\copies{\str C}{\str A}$ by $m$-tuples of colours, while pairing a colouring of $\copies{\str C}{\str A}$ with the index in $\{1,\ldots,m\}$ multiplies the number of forced colours by $m$ (taking $\str B=\str A$ when $t^{\rm cp}_{\mathcal C}(\str A)=1$). In particular, $t_{\mathcal C}(\str A)$ is finite if and only if $t^{\rm cp}_{\mathcal C}(\str A)$ is finite, and for a \emph{rigid} structure $\str A$ (one whose only automorphism is the identity) we have $t_{\mathcal C}(\str A)=t^{\rm cp}_{\mathcal C}(\str A)$.
\leanverifiedat{BANANA/NonPrecompact/FirstOrderEmbeddingCopyDegree.lean}{108}{embedding Ramsey degree is the copy degree times the number of automorphisms}

For languages $L\subseteq L^*$, an \emph{expansion} of an $L$-structure enriches it on the same underlying set by interpretations of the symbols in $L^*\setminus L$, and forgetting those symbols recovers its $L$-\emph{reduct}.
Following the Fra\"iss\'e-age setting of \cite{KPT} (compare \cite[Section~5.2]{HK}, which does not require heredity in general), an \emph{expansion class}
$\mathcal C^*$ of $\mathcal C$ is a hereditary class of finite
$L^*$-structures whose $L$-reducts are precisely the members of
$\mathcal C$; equivalently in the applications below, it is the age of an
expansion of the Fra\"iss\'e limit.
It is \emph{reasonable} if every
embedding $\str A\to \str B$ in $\mathcal C$ lifts, from every prescribed expansion
of $\str A$ in $\mathcal C^*$, to an embedding into some expansion of $\str B$ in $\mathcal C^*$, and has the \emph{expansion property} relative to $\mathcal C$ if for every $\str A\in\mathcal C$ there is $\str B\in\mathcal C$ whose every expansion in $\mathcal C^*$ embeds every expansion of $\str A$ in $\mathcal C^*$.
For $\str B^*\in\mathcal C^*$ with $L$-reduct $\str B$ and a copy
$\str A'\in\copies{\str B}{\str A}$, write $\langle A'\rangle_{\str B^*}$ for the
$L^*$-substructure of $\str B^*$ generated by $A'$; by heredity it belongs to
$\mathcal C^*$. An \emph{expansion type of $\str A$} is the isomorphism type of
such a pair $(\langle A'\rangle_{\str B^*},\str A')$, and $\mathcal C^*$ is
\emph{precompact} if every $\str A\in\mathcal C$ has only finitely many expansion
types. Since all structures are finite, this is equivalent to copies of
$\str A$ generating only finitely many members of $\mathcal C^*$ up to
isomorphism. When $L^*\setminus L$ consists of relation symbols,
$\langle A'\rangle_{\str B^*}$ has underlying set $A'$, and the expansion types of
$\str A$ are simply the isomorphism types in $\mathcal C^*$ with $L$-reduct
$\str A$. If functions or constants are added, a copy of $\str A$ may generate
a larger substructure, and all these substructures are counted.
A \emph{Ramsey expansion} of $\mathcal C$ is a reasonable expansion class with the embedding Ramsey property.
Throughout, \emph{precompact Ramsey expansion} means this hereditary Fra\"iss\'e-age notion; we make the convention explicit because heredity is used in \Cref{lem:expansion-degree-bound}.

\begin{lemma}\label{lem:expansion-degree-bound}
If $\mathcal C$ has a precompact Ramsey expansion and $\str A$ has $t$ expansion
types, then $t^{\rm cp}_{\mathcal C}(\str A)\le t$.
\end{lemma}
\ValidationRelative{\leanrepo/blob/\leancommit/BANANA/NonPrecompact/NativePointedRamseyExpansion.lean\#L72}%
  {native copy-degree bound from a generated pointed cover; many-sorted syntax correspondence}
\ValidationDetails{%
  \leanverifiedat{BANANA/NonPrecompact/FirstOrderGeneratedRamseyExpansion.lean}{112}{copy Ramsey degree bound by the number of pointed generated expansion types in a single-sorted language}%
  \leanverifiedat{BANANA/NonPrecompact/FirstOrderHereditaryExpansionRestriction.lean}{35}{hereditary restriction along every reduct embedding when no function symbols are added}%
  \leanverifiedat{BANANA/NonPrecompact/FirstOrderRamseyExpansion.lean}{290}{infinite copy Ramsey degree excludes precompact expansions by added relations}%
}%
\begin{proof}
Fix $\str B\in\mathcal C$, an expansion $\str B^*\in\mathcal C^*$ of $\str B$,
and a number of colours $r$. If $\str B$ contains no copy of $\str A$,
take $\str C=\str B$. Otherwise, let
$(\str D_1^*,\str A_1'),\ldots,(\str D_k^*,\str A_k')$,
$k\le t$, represent the expansion types $(\langle A'\rangle_{\str B^*},\str A')$
for $\str A'\in\copies{\str B}{\str A}$. Apply the Ramsey property of
$\mathcal C^*$ successively to obtain $\str C_1^*\longrightarrow(\str B^*)^{\str D_1^*}_r$
and $\str C_i^*\longrightarrow(\str C_{i-1}^*)^{\str D_i^*}_r$ for $1<i\le k$, and let
$\str C$ be the $L$-reduct of $\str C_k^*$.
Given an $r$-colouring $\chi$ of $\copies{\str C}{\str A}$, colour each embedding
$\hat e$ of $\str D_i^*$ by $\chi(\hat e(A_i'))$; here $\hat e(A_i')$ is a copy
of $\str A$ in $\str C$, since $\hat e$ is also an embedding of $L$-reducts.
Working backwards through the witnesses gives an embedding
$f\colon\str B^*\to\str C_k^*$ such that, for each $i$, all embeddings
$f\circ e$ with $e\in\Emb(\str D_i^*,\str B^*)$ have the same colour.
Every copy of $\str A$ in $f(B)$ has the form $f(A')$ with
$\str A'\in\copies{\str B}{\str A}$, and an isomorphism
$(\str D_i^*,\str A_i')\to(\langle A'\rangle_{\str B^*},\str A')$ gives
$e\in\Emb(\str D_i^*,\str B^*)$ with $f(e(A_i'))=f(A')$. Hence the colour of
$f(A')$ depends only on $i$, and at most $k\le t$ colours occur in the copy
$f(B)$ of $\str B$.
\end{proof}

\paragraph{EPPA and dynamics.}
A \emph{partial automorphism} of a structure $\str A$ is an isomorphism between two substructures of $\str A$.
A class $\mathcal C$ has the \emph{extension property for partial
automorphisms} (\emph{EPPA}, also called the \emph{Hrushovski property}) if every structure $\str A\in\mathcal C$ has a finite extension
$\str B\in\mathcal C$ in which each partial automorphism of $\str A$
extends to an automorphism of $\str B$.
It has \emph{coherent EPPA} if these extensions can be chosen as a map
$\Theta$ with $\Theta(qp)=\Theta(q)\Theta(p)$ whenever the range of $p$
is exactly the domain of $q$ \cite{HKN22}.

Automorphism groups carry the topology of pointwise convergence.
A topological group is \emph{Polish} if its topology is separable and
completely metrisable. A topological group $G$ is \emph{amenable} if every continuous action of $G$ on a nonempty compact
Hausdorff space (a compact \emph{$G$-flow}) admits a $G$-invariant regular Borel probability measure, and
\emph{extremely amenable} if every compact $G$-flow has a $G$-fixed point.
A subgroup $H\le G$ is \emph{co-precompact} if the quotient right uniformity
on $G/H$ is precompact, i.e., for every identity neighbourhood $U\subseteq G$ there is a finite subset $F\subseteq G$ with $G=UFH$ \cite[Section~1]{MNT}; for permutation groups (taking $U$ to be a tuple stabiliser) this is equivalent to every
$G$-orbit of finite tuples splitting into finitely many $H$-orbits.
\ValidationRelative{\leanrepo/blob/\leancommit/BANANA/NonPrecompact/FirstOrderOrbitCoPrecompact.lean\#L43}%
  {pointwise co-precompactness and finite splitting of tuple orbits for first-order automorphism groups}
A compact $G$-flow is \emph{minimal} if every orbit is dense.
Every topological group $G$ admits a \emph{universal minimal flow} $M(G)$ (unique up to $G$-flow isomorphism), namely a minimal compact $G$-flow that maps
continuously and $G$-equivariantly onto every minimal compact $G$-flow \cite{KPT,PestovDynamics}.
Following Kechris--Rosendal \cite[Section~6]{KR}, $G$ has
\emph{ample generics} if each diagonal conjugation action on $G^n$,
$n\ge1$, has a comeagre orbit. A subset of a topological space is \emph{nowhere dense} if
its closure has empty interior, \emph{meagre} if it is a countable union of nowhere dense sets, and \emph{comeagre} if its complement is meagre.

\paragraph{Elementary notation.}
As customary, set $[n]=\{1,\ldots,n\}$.
Throughout, $\F=\{0,1\}$ denotes the field with two elements. 
For an element $u\in\F$, $[u]\in\{0,1\}$ denotes its integer representative. 
We write
$
 \Fin=\{X\subseteq\mathbb Q:X\text{ is finite}\}.
$
Note that this set has the structure of an $\F$-vector space with symmetric difference as addition. 
For a finite-dimensional $\F$-space $V$, let $V^*=\operatorname{Hom}_{\F}(V,\F)$ be the space of $\F$-linear maps from $V$ to $\F$.
Any vector space isomorphism $h:V\to V$ acts on $V^*$ by
$(h^{-1})^*\phi=\phi\circ h^{-1}$. 
For vector spaces $A$ and~$B$ equipped with bilinear forms, we write $A\perp B$ for the direct sum of $A$ and $B$, equipped with the bilinear form that extends the forms on $A$ and $B$ and sets all cross-values to zero.

\section{Three classes and their Fra\"iss\'e limits}\label{sec:class}

In this section, we introduce the three classes of structures studied in this paper.
\subsection{\BoolBANANA: Boolean algebras with parity}
\label{sec:pre-banana}

We write $\vee,\wedge,\neg,0,1$ for the join, meet, complement, least
element and greatest element of a Boolean algebra.
In a Boolean algebra or ring, $x\le y$ means
$x\wedge y=x$. 
An \emph{atom} is a minimal nonzero element for this order. 
The set of atoms of a Boolean algebra $\str A$ is denoted by $\At(\str A)$.

\begin{definition}\label{def:pre-banana}
A \emph{\BoolBANANA{} structure} is a finite Boolean algebra $\str A$, with
$0\ne1$, equipped with a map $\mu_{\str A}:A\to\F$ satisfying
\[
 \mu_{\str A}(0)=\mu_{\str A}(1)=0,\qquad
 \mu_{\str A}(x\vee y)=\mu_{\str A}(x)+\mu_{\str A}(y)\quad\text{if }x\wedge y=0.
\]
Its language consists of the Boolean operations and constants and a unary
predicate $P(x)$ meaning $\mu_{\str A}(x)=1$. Embeddings are unital
Boolean-algebra embeddings preserving $\mu$. We write $\Kbool$ for this
class.
\end{definition}
\ValidationRelative{\leanrepo/blob/\leancommit/BANANA/NonPrecompact/PrebananaFiniteBooleanCoordinates.lean\#L64}%
  {atom-repr}

We call an atom $a$ \emph{marked} if $\mu_{\str A}(a)=1$. 
For $x \in A$, the additivity yields
\[
 \mu_{\str A}(x)=\sum_{\substack{a\in\At(\str A)\\a\le x}}\mu_{\str A}(a)
 =\bigl|\{a\in\At(\str A):a\le x,\ \mu_{\str A}(a)=1\}\bigr|\pmod2.
\]
\leanverifiedat{BANANA/NonPrecompact/PrebananaAbstractAtomsLiteral.lean}{46}{atom-sum}
In particular, the total number of marked atoms is even.
\leanverifiedat{BANANA/NonPrecompact/PrebananaAbstractAtomsLiteral.lean}{61}{even-atoms} If $f \colon \str A\to \str B$ is an
embedding of Boolean algebras, set
$V_a=\{b\in\At(\str B):b\le f(a)\}
$ for all $a \in \At(\str A)$.
The sets $V_a$ form a partition of $\At(\str B)$ into nonempty blocks
\leanverifiedat{BANANA/NonPrecompact/PrebananaAbstractMarkedPartitions.lean}{125}{atom-fibres}
\ValidationDetails{%
\leanverifiedat{BANANA/NonPrecompact/PrebananaFibrewiseMarkedAllocation.lean}{27}{base-compatible fibre marks}
}%
and
$\sum_{b\in V_a}\mu_{\str B}(b)=\mu_{\str A}(a)$. Conversely, any such partition defines
an embedding $f \colon \str A \to \str B$ by setting
$f(x)=\bigvee\{b:b\in V_a\text{ for some atom }a\le x\}$ for all $x \in A$.

\begin{proposition}\label{prop:pre-banana-fraisse}
The class $\Kbool$ is a strong amalgamation Fra\"iss\'e class.
Its Fra\"iss\'e limit $\Mbool$ is $\omega$-categorical.
\end{proposition}
\leanverifiedat{BANANA/NonPrecompact/BooleanBananaFirstOrderAge.lean}{252}{the marked Boolean cylinder is the Fra\"iss\'e limit of the actual marked Boolean age}
\ValidationDetails{%
  \leanpartialat{BANANA/NonPrecompact/PrebananaFiveAryCylinderFiniteAge.lean}{43}{finite age embeddings}%
}%
\begin{proof}
Let $\str A,\str B, \str C \in \Kbool$ be Boolean algebras with embeddings of $\str A$ into $\str B$ and $\str C$. 
We identify $\str A$ with a common substructure of $\str B$ and $\str C$.

For every $a \in \At(\str A)$, let $b_1^a, \dots, b_{s_a}^a$ and
$c_1^a, \dots, c_{t_a}^a$ be the atoms of $\str B$ and $\str C$ below $a$, respectively.
Let $\str D$ be a Boolean algebra with atoms
$\{d^a_{ij} \colon a\in\At(\str A),\ (i,j) \in [s_a] \times [t_a]\}$.
For a fixed $a\in\At(\str A)$, we regard $(d^a_{ij})_{i,j}$ as an
$(s_a \times t_a)$-matrix.
Let
$r_i^a=\mu_{\str B}(b_i^a)$ and $u_j^a=\mu_{\str C}(c_j^a)$ in $\F$, and define
\[
 \mu_{\str D}(d^a_{ij})=
 \begin{cases}
  0,&i<s_a,\ j<t_a,\\
  r_i^a,&i<s_a,\ j=t_a,\\
  u_j^a,&i=s_a,\ j<t_a,\\
  r_{s_a}^a+\sum_{j<t_a}u_j^a,&i=s_a,\ j=t_a
 \end{cases},
\]
and extend $\mu_{\str D}$ additively to all elements of $D$.
\leanpartial{BANANA/NonPrecompact/PrebananaIndexedParity.lean}%
  {finite-grid}
For all $i \in [s_a]$, we have
$\sum_{j = 1}^{t_a} \mu_{\str D}(d^a_{ij}) = r_i^a$. 
For all $j \in [t_a]$, we have 
$\sum_{i= 1}^{s_a} \mu_{\str D}(d^a_{ij}) = u_j^a$. 
To see this for $ j = t_a$, note that
$\sum_{i = 1}^{s_a} r_i^a=\mu_{\str A}(a)=\sum_{j=1}^{t_a} u_j^a$.
\ValidationDetails{%
\leanverifiedat{BANANA/NonPrecompact/ParityStarMatrix.lean}{37}{sum_row_parityStarMatrix}
\leanverifiedat{BANANA/NonPrecompact/ParityStarMatrix.lean}{82}{sum_col_parityStarMatrix}
\leanverifiedat{BANANA/NonPrecompact/FiberProductParityAmalgam.lean}{39}{sum_fiberGridMark_right}
\leanverifiedat{BANANA/NonPrecompact/FiberProductParityAmalgam.lean}{62}{sum_fiberGridMark_left}
}%

Define maps $f_B \colon B \to D$ and $f_C \colon C \to D$ by setting
$f_B(b_i^a) = \bigvee_j d^a_{ij}$ and $f_C(c_j^a) = \bigvee_i d^a_{ij}$ for every $a\in\At(\str A)$, and extending
these assignments additively. The above computation of the row and column sums of $\mu_{\str D}$ shows that $f_B$ and $f_C$ preserve $\mu$.
\ValidationDetails{%
\leanverifiedat{BANANA/NonPrecompact/PrebananaGridMarkedEmbeddings.lean}{134}{grid-left}
\leanverifiedat{BANANA/NonPrecompact/PrebananaGridMarkedEmbeddings.lean}{147}{grid-right}
}%
\leanverifiedat{BANANA/NonPrecompact/PrebananaPowersetBridge.lean}{117}{all-elements}
For every $a\in\At(\str A)$, a union of rows of the matrix $(d_{ij}^a)$ can only equal a
union of columns if both are empty or both are the full matrix.
Thus $f_B(B) \cap f_C(C) = f_B(A) = f_C(A)$, which proves strong amalgamation. 
\ValidationDetails{%
\leanverifiedat{BANANA/NonPrecompact/FiberGridCoordinates.lean}{194}{intersection}
\leanverified{BANANA/NonPrecompact/PrebananaAbstractStrongAmalgamation.lean}{abstract AP}
}%
\leanverified{BANANA/NonPrecompact/PrebananaFiniteBooleanModelStrongAP.lean}{finite model strong AP}

Heredity and countability of $\Kbool$ are immediate. 
\leanverified{BANANA/NonPrecompact/PrebananaFiniteBooleanModelHeredity.lean}{finite-model heredity}
\ValidationDetails{%
\leanpartial{BANANA/NonPrecompact/PrebananaFiniteClass.lean}%
  {countability}
}%
\leanverified{BANANA/NonPrecompact/PrebananaFiniteBooleanModelSkeleton.lean}{canonical finite skeleton}
The joint embedding property is obtained by amalgamation over
the Boolean algebra $\{0,1\}$ whose parity function is the constant zero map.
\ValidationDetails{%
\leanpartial{BANANA/NonPrecompact/PrebananaFiniteJEP.lean}%
  {atom-JEP}
}%
\leanverified{BANANA/NonPrecompact/PrebananaFiniteBooleanModelJEP.lean}{finite-model JEP}
Moreover, $\Kbool$ contains arbitrarily large structures.
\ValidationDetails{%
\leanpartial{BANANA/NonPrecompact/PrebananaFiniteClass.lean}%
  {unboundedness}
}%
An $n$-tuple of elements of $\Mboolset$ generates a Boolean algebra with at most
$2^n$ atoms as an atom is specified by membership or nonmembership in each
of the $n$ generators.
\ValidationDetails{%
\leanverified{BANANA/NonPrecompact/PrebananaGeneratedBooleanClosure.lean}{term-value closure}
\leanpartial{BANANA/NonPrecompact/PrebananaTupleSignatures.lean}%
  {pattern-bound}
}%
\leanverified{BANANA/NonPrecompact/PrebananaAbstractMarkedTupleProfile.lean}{abstract profiles}
Only finitely many generated algebras, parity maps
and distinguished tuples are possible. Homogeneity gives finitely many
orbits on $n$-tuples, proving $\omega$-categoricity.
\ValidationDetails{%
\leanverified{BANANA/NonPrecompact/PrebananaAbstractGeneratedMarkedIso.lean}{abstract iso}
\leanverified{BANANA/NonPrecompact/PrebananaAbstractGeneratedBooleanOperations.lean}{Boolean operations}
\leanverified{BANANA/NonPrecompact/PrebananaActualMarkedEmbedding.lean}{marked Boolean embedding}
\leanverifiedat{BANANA/NonPrecompact/PrebananaFiveAryCylinderUltrahomogeneity.lean}{49}{global homogeneity}
\leanverifiedat{BANANA/NonPrecompact/PrebananaFiveAryCylinderMarkedOrbitFiniteness.lean}{91}{finite tuple orbits}%
\leanverifiedat{BANANA/NonPrecompact/BooleanBananaFirstOrderCategoricity.lean}{104}{the marked Boolean cylinder limit is countably infinite}%
  \leanverifiedat{BANANA/NonPrecompact/BooleanBananaFirstOrderAge.lean}{201}{strong amalgamation of the actual marked Boolean age}%
  \leanverifiedat{BANANA/NonPrecompact/BooleanBananaFirstOrderAge.lean}{252}{the marked Boolean cylinder is the limit of its native finite class}%
}%
\leanverifiedat{BANANA/NonPrecompact/BooleanBananaFirstOrderCategoricity.lean}{96}{countable categoricity of the marked Boolean cylinder limit}
\end{proof}

\subsection{\FolkBANANA: disjoint-union embeddings}
\label{sec:folkman-class}

We now define a variant of the previous construction that considers Boolean rings instead of Boolean algebras. 
The resulting geometry underlies Folkman's theorem \cite{GRS}.
For this reason, we call the resulting structures \FolkBANANA{} structures.

A finite \emph{Boolean ring} is a finite Boolean algebra in the language
$(\triangle,\wedge,0)$ without a named greatest element. 
We also allow the
one-element zero ring. 
For the ring structure, the addition is given by symmetric difference, so
$x+y=x\triangle y$, and the multiplication is meet, so $xy=x\wedge y$.

An embedding $f \colon \str A \to \str B$ between Boolean rings $\str A$ and $\str B$ is an injective ring homomorphism, i.e., an injective map compatible with $\triangle$ and $\wedge$. 
A ring homomorphism is not necessarily unital, so in particular, an embedding does not need to preserve the greatest element.
Furthermore, note that joins are defined by addition and multiplication by
$x\vee y=(x\triangle y)\triangle(x\wedge y)$.

\begin{definition}\label{def:folkman-banana}
A \emph{\FolkBANANA{} structure} consists of finite Boolean rings
$L_A,R_A$ on disjoint sorts and a map $\beta_{\str A}:L_A\times R_A\to\F$ that is
additive in each variable.
 We encode $\beta_{\str A}$ by a binary relation 
$E_{\str A}(x,y)\Longleftrightarrow\beta_{\str A}(x,y)=1$.
An embedding of $\str A$ into $\str B$ is a pair $f = (f_L, f_R)$, where $f_L \colon L_A \to L_B$ and $f_R \colon R_A \to R_B$ are embeddings of Boolean rings, such that~$f$ preserves $\beta$.
Write $\Kfolk$ for the class of \FolkBANANA{} structures and their isomorphic copies.
\end{definition}
\ValidationDetails{%
\leanverified{BANANA/NonPrecompact/FolkmanNativeNonUnitalRingOwners.lean}{native non-unital owners}
}%

For Boolean rings $\str A,\str B$ and an embedding $f:\str A\to \str B$, the sets
$
 V_a=\{b\in\At(\str B):b\le f(a)\}
$
for $a \in \At(\str A)$
are pairwise disjoint nonempty blocks. 
\ValidationDetails{%
\leanverified{BANANA/NonPrecompact/FolkmanRingPartialOwners.lean}{partial atom owners}
}%
Conversely, for all $x \in A$, we have $f(x) = \bigvee\{b\in V_a:a\in\At(\str A),\ a\le x\}$, so these blocks determine $f$.
\leaninterfaceat{BANANA/NonPrecompact/FolkmanNativeFiniteClassMorphismWitness.lean}%
  {25}{atom codes for the given maps}
In contrast to the analogue for \BoolBANANA{} structures, the blocks do not need to cover $\At(\str B)$, i.e., we do not necessarily have $\At(\str B) \subseteq \bigcup_{a \in \At(\str A)} V_a$.
This, rather than the passage to two sorts, is the essential difference
between the two classes; a single Boolean ring with an additive map to $\F$
would serve equally well. We use two sorts because they make the Ramsey
obstacle more explicit, and because forgetting the meet operation gives two
$\F$-vector spaces with a bilinear map, that is, a BANANA structure.

\begin{proposition}\label{prop:folkman-fraisse}
The class $\Kfolk$ is a strong amalgamation Fra\"iss\'e class.
\end{proposition}
\begin{proof}
Let $\str A,\str B,\str C \in \Kfolk$ be such that $\str A$ embeds into $\str B$ and $\str C$. 
We identify $\str A$ with a common substructure of $\str B$ and $\str C$. 
First we amalgamate the left ring sort.
For each $a\in\At(L_A)$, let 
$ (L_B)_a=\{b\in\At(L_B):b\le a\}$ and $
 (L_C)_a=\{c\in\At(L_C):c\le a\}$, and set $
 D_a=(L_B)_a\times (L_C)_a$.
Put
\[
 \Omega_L=
 \biggl(\coprod_{a\in\At(L_A)}D_a\biggr)
 \sqcup
 \{\hat b:b\in\At(L_B),\ b\not\le1_{L_A}\}
 \sqcup
 \{\hat c:c\in\At(L_C),\ c\not\le1_{L_A}\},
\]
with all introduced points being distinct. 
Let $L_D=\mathcal P(\Omega_L)$.
Define embeddings
$i_B^L\colon L_B\to L_D$ and $i_C^L\colon L_C\to L_D$ as follows.
For all $a \in \At(L_A)$, we map $b\in (L_B)_a$ to
$\{b\}\times (L_C)_a$ and $c\in (L_C)_a$ to $(L_B)_a\times\{c\}$. 
For all $b \in \At(L_B)$ with $b \not \leq 1_{L_A}$, we map $b$ to $\{\hat{b}\}$, and similarly for $c \in \At(L_C)$ with $c \not \leq 1_{L_A}$.
Arguing similarly to the proof of
\Cref{prop:pre-banana-fraisse} shows that the images of $L_B$ and $L_C$ in $L_D$ intersect exactly in
$L_A$. 
\ValidationDetails{%
\leanverified{BANANA/NonPrecompact/FolkmanAtomAmalgam.lean}{atom grid}
}%
\leanverified{BANANA/NonPrecompact/FolkmanNativeSortwiseIntersection.lean}{sortwise AP}
Applying this construction again for the right sort yields a Boolean ring $R_D$ with
embeddings $i_B^R\colon R_B\to R_D$ and $i_C^R\colon R_C\to R_D$.
\ValidationDetails{%
\leanverified{BANANA/NonPrecompact/FolkmanFiniteBooleanRingMaps.lean}{ring maps, intersection}%
\leanverified{BANANA/NonPrecompact/FolkmanFiniteTwoSortRingDiagram.lean}{two-sort ring diagram}
\leanverified{BANANA/NonPrecompact/FolkmanNativeTwoSortOwnerDiagram.lean}{four native owner maps}
\leanverified{BANANA/NonPrecompact/FolkmanIntersectionTransport.lean}{native carrier intersection}
\leanverified{BANANA/NonPrecompact/FolkmanNativeIntersection.lean}{exact original-ring intersection}
\leanverified{BANANA/NonPrecompact/FolkmanNativeRingAtomEquiv.lean}{native ring atom coordinates}
\leanverified{BANANA/NonPrecompact/FolkmanNativeRingAtomInverse.lean}{inverse ring coordinates}
\leanverified{BANANA/NonPrecompact/FolkmanNativeRingMapCoordinates.lean}{non-unital coordinate naturality}
\leanverified{BANANA/NonPrecompact/FolkmanNativeEmbeddings.lean}{native ring embedding operations}
\leanverified{BANANA/NonPrecompact/FolkmanAtomPivots.lean}{actual pivot retractions}
\leanverified{BANANA/NonPrecompact/FolkmanPivotCrossIdentities.lean}{pivot cross identities}
}%

We identify $L_A$ with its images in $L_B$ and $L_C$, and similarly for $R_A$.
Choose $\F$-vector space
complements $L_B=L_A\oplus U_B$, $L_C=L_A\oplus U_C$ and
$R_B=R_A\oplus V_B$, $R_C=R_A\oplus V_C$. 
Retain $\beta_{\str B}$ on $L_B\times R_B$ and $\beta_{\str C}$ on $L_C\times R_C$,
and set
$\beta(U_B,V_C)=\beta(U_C,V_B)=0$. 
The prescriptions agree on the common base,
so they define a bilinear map. Extend it to the full amalgam rings $L_D$ and $R_D$ by
choosing vector space complements and putting zero on the remaining terms.
\ValidationDetails{%
\leanverified{BANANA/NonPrecompact/FolkmanPairingExtensionCertificate.lean}{split pairing certificate}
\leanverified{BANANA/NonPrecompact/FolkmanActualTwoSortPairing.lean}{actual-grid pairing}
\leanverified{BANANA/NonPrecompact/FolkmanNativePairingTransport.lean}{native pairing transport}
}%
\leanverifiedat{BANANA/NonPrecompact/FolkmanNativePairingCompatibility.lean}{27}{pairing compatibility}
This is a strong amalgam. 
Heredity and countability are immediate. 
Amalgamating over the two zero rings yields the joint embedding property. 
Clearly, $\Kfolk$ also contains arbitrarily large structures.
This shows that $\Kfolk$ is a strong amalgamation Fra\"iss\'e class.
\leanverifiedat{BANANA/NonPrecompact/FolkmanNativeCompleteFiniteClass.lean}%
  {85}{native finite class}
\ValidationDetails{%
\ValidationRelative{\leanrepo/blob/\leancommit/BANANA/NonPrecompact/FolkmanFirstOrderFraisse.lean\#L146}%
  {first-order Fra\"iss\'e class in the product encoding}
  \leanverifiedat{BANANA/NonPrecompact/FolkmanFirstOrderFraisse.lean}{93}{strong amalgamation with exact intersections in the first-order product encoding}%
}%
\ValidationDetails{%
\leanverifiedat{BANANA/NonPrecompact/FolkmanNativeFiniteHeredity.lean}{145}{hereditary restriction of the native non-unital pairings}
\leanverified{BANANA/NonPrecompact/FolkmanFiniteClassData.lean}{finite class data}
\leanverified{BANANA/NonPrecompact/FolkmanNativeJEP.lean}{joint embedding}
\leanverified{BANANA/NonPrecompact/FolkmanFiniteStrongAmalgamPairing.lean}{finite-atom strong AP}
\leanverified{BANANA/NonPrecompact/FolkmanNativeStrongAmalgam.lean}{native six-ring pairing certificate}
\leanverified{BANANA/NonPrecompact/FolkmanNativeStrongAmalgamPackage.lean}{native strong amalgamation}
\leanverified{BANANA/NonPrecompact/FolkmanNativePairingAtomMatrix.lean}{native atom-matrix reconstruction}
}%
\end{proof}

\begin{proposition}\label{prop:folkman-omega-categorical}
The Fra\"iss\'e limit $\Mfolk$ of $\Kfolk$ is $\omega$-categorical.
\end{proposition}
\begin{proof}
A Boolean ring generated by $x_1,\ldots,x_n$ has at most $2^n-1$ atoms.
\leanpartialat{BANANA/NonPrecompact/FolkmanTupleProfiles.lean}{67}{nonzero atom patterns}
\ValidationDetails{%
\leanverifiedat{BANANA/NonPrecompact/FolkmanFirstOrderOrbits.lean}{99}{sharp atom bound for generated finite representatives}%
}%
Indeed put $u=x_1\vee\cdots\vee x_n$. For each
$\varnothing\ne S\subseteq[n]$, form
\[
 c_S=\bigwedge_{i\in S}x_i\ \wedge\!
       \bigwedge_{j\notin S}(u\triangle x_j).
\]
The nonzero $c_S$ are exactly the atoms of the generated ring; the omitted
all-zero membership pattern lies outside $u$.
\ValidationDetails{%
\leanverified{BANANA/NonPrecompact/FolkmanTupleProfilesPairingTransport.lean}{tuple profiles}
}%
An $n$-tuple from $\Mfolkset$ therefore generates a structure of bounded
size. As in \Cref{prop:pre-banana-fraisse}, there are only finitely many
such structures with a distinguished tuple and a map between the two sorts.
Homogeneity gives finitely many tuple orbits.
\ValidationRelative{\leanrepo/blob/\leancommit/BANANA/NonPrecompact/FolkmanFirstOrderExistence.lean\#L24}%
  {existence and countable categoricity of the encoded Folkman--BANANA limit}
\ValidationDetails{%
\leanverifiedat{BANANA/NonPrecompact/FolkmanFirstOrderCategoricity.lean}{52}{the supplied encoded limit is countably infinite}%
\ValidationRelative{\leanrepo/blob/\leancommit/BANANA/NonPrecompact/FolkmanFirstOrderTaggedLimitOrbits.lean\#L15}%
  {finite native mixed tuple orbits of the encoded Fra\"iss\'e limit}%
\leanverifiedat{BANANA/NonPrecompact/FolkmanFirstOrderDecodingLimit.lean}{14}{native limit factorisation}%
\leanverifiedat{BANANA/NonPrecompact/FolkmanFirstOrderOrbits.lean}{171}{finite encoded tuple orbit quotient}%
}%
\end{proof}

\subsection{BANANA: the bilinear class}
\label{sec:banana-class}

\begin{definition}\label{def:banana}
A \emph{BANANA structure} is a triple
$\str A=(L_A,R_A,\beta_{\str A})$, where $L_A$ and $R_A$ are finite-dimensional
$\F$-vector spaces on disjoint sorts and $
 \beta_{\str A}:L_A\times R_A\longrightarrow\F$
is bilinear, possibly degenerate.
Its language has addition and zero on
each sort and a binary relation
$E_{\str A}(x,y)\Longleftrightarrow\beta_{\str A}(x,y)=1$. An embedding is a pair of
injective linear maps preserving $\beta_{\str A}$ or, equivalently, $E_{\str A}$.
Write $\K$ for this class.
\end{definition}
\leaninterfaceat{BANANA/NonPrecompact/BananaStructure.lean}%
  {71}{BananaMatrixEmbedding.edge_iff}

Let $\str A = (L_A, R_A, \beta_{\str A})$ be a BANANA structure.
Then $|A|=2^{\dim L_A}+2^{\dim R_A}$. 
A substructure of $\str A$ is a structure $\str A'=(L_{A'}, R_{A'}, \beta_{\str A'})$, where $L_{A'}$ and $R_{A'}$ are subspaces of $L_A$ and $R_A$, respectively, and $\beta_{\str A'} = \beta_{\str A}|_{L_{A'} \times R_{A'}}$. 
The structure $\str A$ is called a \emph{bilinear pairing} if there exist bases $\{l_1, \dots, l_{\dim L_A}\}$ and $\{r_1, \dots, r_{\dim R_A}\}$ of $L_A$ and $R_A$, respectively, such that $\beta_{\str A}(l_i, r_j) = 1$ only if $i = j$. 
We call a pairing \emph{perfect} if also the converse holds, i.e., if $\dim L_A = \dim R_A$ and $\beta_{\str A}(l_i, r_j) = 1$ exactly if $i = j$.
Intuitively, if $\str A$ is a pairing, then~$E_{\str A}$ defines a matching between $\{l_1, \dots, l_{\dim L_A}\}$ and $\{r_1, \dots, r_{\dim R_A}\}$, and if $\str A$ is a perfect pairing, then~$E_{\str A}$ defines a perfect matching.
The property of being a perfect pairing can easily be phrased in terms of dual vector spaces. 
Clearly, there is a natural map $\Phi_L \colon L_A \to R_A^*$, mapping $l \in L_A$ to $\varphi_l \colon R_A \to \F,\, r \mapsto \beta_{\str A}(l,r)$.
Similarly, we define a natural map $\Phi_R \colon R_A \to L_A^*$. 
Then $\str A$ is a perfect pairing if and only if the maps $\Phi_L$ and $\Phi_R$ are isomorphisms.

\begin{lemma}\label{lem:completion}
Every BANANA structure $\str A=(L_A,R_A,\beta_{\str A})$ embeds in a perfect
pairing with both sides of dimension $\dim L_A+\dim R_A$.
\end{lemma}
\leanverifiedat{BANANA/NonPrecompact/BananaNativePerfectPairing.lean}{62}{native perfect completion using direct sums with the dual spaces}
\ValidationDetails{%
  \leaninterfaceat{BANANA/NonPrecompact/Completion.lean}%
    {80}{BananaMatrixStructure.completionEmbedding}%
}%
\begin{proof}
Let $\widehat L=L_A\oplus R_A^*$, $\widehat R=R_A\oplus L_A^*$ and
define $
 \widehat\beta((x,\alpha),(y,\phi))
 =\beta_{\str A}(x,y)+\alpha(y)+\phi(x)$ for all $(x,\alpha) \in \widehat L$ and $(y,\phi) \in \widehat R$.
It is easily verified that $(\widehat L, \widehat R, \widehat{\beta})$ is a perfect pairing.
The structure $\str A$ embeds into $(\widehat{L}, \widehat{R}, \widehat{\beta})$ by mapping $x \in L_A$ to $(x,0) \in \widehat L$ and $y \in R_A$ to $(y,0) \in \widehat R$.
\end{proof}

For every finite set $S\subseteq\mathbb Q$, set
$ \Bban{S}=(\F^S,\F^S,\beta_S)$, where $
 \beta_S(x,y)=\sum_{s\in S}x_s y_s$ for $x,y \in \F^S$.
We also write $\Bban{D} = \Bban{[D]}$. 
If $\str A = (L_A, R_A, \beta_{\str A})$ is a perfect pairing with $\dim L_A = \dim R_A = D$ for some $D \in \mathbb{N}$, choosing a basis of $L_A$ and the corresponding dual basis for $R_A$ identifies~$\str A$ with $\Bban{D}$. 
If we identify a vector of $\F^S$ with its support in $S$, then $\beta_S(x,y)$ is the parity of the intersection of the supports of
$x$ and $y$.
Together with~\Cref{lem:completion}, this implies that every element of $\K$ embeds into
\begin{equation*}
 \str M=(L,R,+_L,+_R,0_L,0_R,E),
\end{equation*}
where $L=R=\Fin$, we set $X +_L Y := X \triangle Y$ for $X, Y \in L$ (and similarly for $R$), and 
$E(X,Y)\Longleftrightarrow|X\cap Y|\text{ is odd}$.

We use the following elementary fact.

\begin{lemma}\label{lem:extendlinearmap}
Let $V$ be a finite-dimensional $\F$-vector space, let $q,q':V\to\F^k$ be surjective
linear maps, and let $f:A\to A'$ be an isomorphism between subspaces $A, A'$ of $V$,
with $q'f=q$ on $A$. Then $f$ extends to $h\in\mathrm{GL}(V)$ with
$q'h=q$ on $V$.
\end{lemma}
\leanverifiedat{BANANA/NonPrecompact/LinearExtension.lean}{62}{exists_linearEquiv_extends_of_surjective}

\begin{lemma}\label{lem:perfect-homogeneous}
Every isomorphism between substructures of a finite perfect pairing
$\str A = (L_A,R_A,\beta_{\str A})$ extends to an automorphism of that pairing.
\end{lemma}
\leanverifiedat{BANANA/NonPrecompact/BananaNativePerfectPairing.lean}{89}{extension of every subpair isomorphism in an arbitrary finite perfect pairing}
\ValidationDetails{%
  \ValidationRelative{\leanrepo/blob/\leancommit/BANANA/NonPrecompact/PerfectHomogeneity.lean\#L157}%
  {exists_perfectPairAutomorphism_extends}%
}%
\begin{proof}
As before, we identify $R_A$ with $L_A^*$.
Let $\str U = (L_U, R_U, \beta_{\str U})$ and $\str U' = (L_{U'}, R_{U'}, \beta_{\str U'})$ be substructures of $\str A$, and let $f \colon \str U \to \str U'$ be an isomorphism. 
It defines vector space isomorphisms $f_L \colon L_U \to L_{U'}$ and $f_R \colon R_U \to R_{U'}$. 
For a basis $b_1,\ldots,b_k$ of $R_U$, set $b_i'=f_R(b_i)$ for $i = 1, \dots, k$. 
For $v \in L_A$, define
$q(v)=(b_i(v))_{i=1}^k$ and $q'(v)=(b_i'(v))_{i=1}^k$.
The maps $q, q' \colon L_A \to \F^k$ are surjective because $b_1, \dots, b_k$ and $b_1', \dots, b_k'$ are linearly independent in $R_A \cong L_A^*$.
Since $E_{\str A}$ needs to be preserved, we have $q'f_L=q$ on $L_U$.
By \Cref{lem:extendlinearmap}, there exists $h \in \operatorname{GL}(L_A)$ extending $f_L$ with $q'h=q$.
Then $(h^{-1})^*b_i=b_i'=f_R(b_i)$, so $(h,(h^{-1})^*)$ extends $f$.
\end{proof}

\begin{proposition}\label{prop:fraisse}
The class $\K$ is a strong amalgamation class with Fra\"iss\'e limit~$\str M$, and~$\str M$ is $\omega$-categorical.
\end{proposition}
\begin{proof}
Let $\str A, \str B, \str C \in \K$ with embeddings of $\str A$ into $\str B$ and $\str C$, and identify $\str A$ with substructures of $\str B$ and $\str C$.
Write 
$L_{B}=L_A\oplus U_B$ and $R_{B}=R_A\oplus V_B$ for vector spaces $U_B$ and $V_B$, and similarly for $\str C$.
Set $L=L_A\oplus U_B \oplus U_C$ and $R=R_A\oplus V_B\oplus V_C$. 
Let $\iota_B^L \colon L_B \to L$ and $\iota_B^R \colon R_B \to R$ be the natural embeddings, and similarly for $\str C$.
Define a bilinear map $\beta \colon L \times R \to \F$ by $\beta|_{\iota_B^L(L_B) \times \iota_B^R(R_B)} = \beta_{\str B}$, $\beta|_{\iota_C^L(L_C) \times\iota_C^R(R_C)} = \beta_{\str C}$ and $\beta(U_B,V_C)=\beta(U_C,V_B)=0$, extended bilinearly.
Then $(L,R,\beta)$ is a strong amalgam of $\str B$ and $\str C$ over $\str A$.
\ValidationDetails{%
\leanpartialat{BANANA/NonPrecompact/BananaAmalgam.lean}%
  {191}{bananaBlockAmalgam_left}
\leanpartialat{BANANA/NonPrecompact/BananaAmalgam.lean}%
  {203}{bananaBlockAmalgam_right}
\leanpartialat{BANANA/NonPrecompact/BananaAmalgam.lean}%
  {132}{sumAmalgam_eq_iff_common}
\leanverifiedat{BANANA/NonPrecompact/BananaStrongBlockAmalgam.lean}%
  {108}{split B pairs}
\leanverifiedat{BANANA/NonPrecompact/BananaStrongBlockAmalgam.lean}%
  {116}{split C pairs}
\leaninterfaceat{BANANA/NonPrecompact/BananaTriplePairing.lean}%
  {218}{strong amalgam in coordinates}%
}%
\leanverifiedat{BANANA/NonPrecompact/BananaNativeStrongAmalgamation.lean}{136}{native strong amalgamation with exact image intersections and vanishing cross blocks}
The joint embedding property is obtained by amalgamating over the zero spaces. 
\leanverifiedat{BANANA/NonPrecompact/BananaNativeStrongAmalgamation.lean}{208}{native joint embedding over the zero pairing}
\ValidationDetails{%
  \leaninterfaceat{BANANA/NonPrecompact/BananaFiniteFraisseProperties.lean}%
    {86}{joint embedding in coordinates}%
}%
That $\K$ is hereditary, countable, and contains structures of unbounded size is clear.
\leanverifiedat{BANANA/NonPrecompact/BananaFiniteFraisseProperties.lean}%
  {27}{hereditary subspaces}
\leanverifiedat{BANANA/NonPrecompact/BananaFiniteCountability.lean}%
  {44}{finite codes}
\leanverifiedat{BANANA/NonPrecompact/BananaFiniteCountability.lean}%
  {52}{unbounded size}

We have shown that $\Age(\str M)=\K$.
\leanverifiedat{BANANA/NonPrecompact/BananaLimitFiniteAge.lean}%
  {108}{finite age}
Let $\str A$ and $\str B$ be finite substructures of $\str M$, and let $f \colon \str A \to \str B$ be an isomorphism. 
Let $S\subseteq\mathbb Q$ be a finite set containing the supports of all vectors in $A$ and $B$. 
Both $\str A$ and $\str B$ lie in the substructure $\Bban{S}$ of $\str M$. 
By \Cref{lem:perfect-homogeneous}, the isomorphism $f$
extends to an automorphism of $\Bban{S}$. 
Extending this automorphism by the identity map on the coordinates not in $S$ yields an automorphism of $\str M$. 
Thus $\str M$ is homogeneous. 
\leanverifiedat{BANANA/NonPrecompact/BananaLimitGlobalHomogeneity.lean}%
  {31}{homogeneity}
\ValidationRelative{\leanrepo/blob/\leancommit/BANANA/NonPrecompact/BananaFunctionalLimit.lean\#L204}%
  {the explicit BANANA model is the Fra\"iss\'e limit in the single-sorted functional encoding}
\ValidationDetails{%
  \leanverifiedat{BANANA/NonPrecompact/BananaFunctionalStrongAge.lean}{39}{strong amalgamation for the functional age of the BANANA limit}%
}%
As any finite set of vectors in $M$ generates a finite $\F$-subspace, $\str M$ is locally finite. 
\leanverifiedat{BANANA/NonPrecompact/BananaLimitLocalFiniteness.lean}%
  {47}{finite spans}

For $\omega$-categoricity, note that an $n$-tuple of elements of $M$ generates a
substructure whose two sides have dimension at most $n$. Up to isomorphism,
there are only finitely many such structures with a distinguished $n$-tuple
(this amounts to the sort of each entry, the linear relations among the entries
of each sort, and the pairings between left and right entries).
By homogeneity, two tuples generating isomorphic pointed structures lie in
the same orbit.
\leanverifiedat{BANANA/NonPrecompact/BananaLimitMixedTupleOrbits.lean}%
  {93}{tagged tuple signatures}
Hence $\Aut(\str M)$ has finitely many orbits on $n$-tuples for every $n$, and
$\str M$ is $\omega$-categorical.
\leanpartialat{BANANA/NonPrecompact/BananaLimitMixedTupleOrbits.lean}%
  {139}{finite tagged tuple orbits}
\ValidationRelative{\leanrepo/blob/\leancommit/BANANA/NonPrecompact/BananaFunctionalLimit.lean\#L187}%
  {countable categoricity in the single-sorted functional language}%
\ValidationDetails{%
  \ValidationRelative{\leanrepo/blob/\leancommit/BANANA/NonPrecompact/BananaLimitFirstOrder.lean\#L138}%
  {countable categoricity of the tagged relational BANANA limit}%
\leanverifiedat{BANANA/NonPrecompact/BananaLimitFirstOrder.lean}{129}{the tagged relational presentation is countably infinite}%
  \ValidationRelative{\leanrepo/blob/\leancommit/BANANA/NonPrecompact/BananaTaggedModelDecoding.lean\#L466}%
  {every model of the tagged theory is isomorphic to the re-encoding of its two vector spaces and bilinear pairing}%
  \leanverifiedat{BANANA/NonPrecompact/BananaSortedSyntax.lean}{302}{translation of two-sorted formulas preserves truth in every decoded tagged BANANA model}%
  \leanverifiedat{BANANA/NonPrecompact/BananaTaggedSyntaxReverse.lean}{403}{translation of tagged formulas back to the two-sorted language preserves truth in every decoded model}%
}%
\end{proof}

\section{EPPA and amenability}\label{sec:eppa}
\subsection{\BoolBANANA}

\ValidationDetails{%
\leanpartialat{BANANA/NonPrecompact/PrebananaEppaObstruction.lean}%
  {21}{no_perm_maps_finset_to_disjoint_union}
}%
\begin{proposition}\label{prop:pre-banana-limits}
The class $\Kbool$ does not have EPPA, and $\Gbool=\Aut(\Mbool)$ is not amenable.
\end{proposition}
\begin{proof}
For EPPA, we consider a Boolean algebra with three unmarked atoms $a,b,c$.
The map sending $a$ to $a\vee b$ extends to a marked isomorphism from
$\{0,a,b\vee c,1\}$ to $\{0,a\vee b,c,1\}$.
No automorphism of a finite
extension can do this: $a\vee b$ has strictly more atoms below it than $a$,
and automorphisms preserve that number.
\ValidationDetails{%
\leanverified{BANANA/NonPrecompact/PrebananaClassEPPA.lean}{atom EPPA}
}%
\leanverifiedat{BANANA/NonPrecompact/PrebananaAbstractClassEPPA.lean}{63}{failure of EPPA in the abstract finite marked Boolean class}

For amenability, use the action of $\Gbool=\Aut(\Mbool)$ on the Stone space
of the underlying Boolean algebra, namely its (compact) space of ultrafilters. Each Boolean element $x$
corresponds to the clopen set $[x]$ of ultrafilters containing it, and $g\cdot[x]=[g(x)]$.
\ValidationDetails{%
\leanverifiedat{BANANA/NonPrecompact/BooleanBananaStoneAction.lean}{146}{continuous Stone action in the pointwise topology}
}%
By homogeneity all proper nonzero elements with $\mu_{\Mbool}(x)=0$ are in the same $\Gbool$-orbit.
\ValidationDetails{%
\leanverifiedat{BANANA/NonPrecompact/PrebananaFiveAryCylinderUltrahomogeneity.lean}{49}{marked homogeneity}
\leanverified{BANANA/NonPrecompact/PrebananaOneGeneratorMarkedIso.lean}{conditional zero-mark orbit}
}%
An invariant probability measure would give their clopen sets a common mass
$t$. The limit contains partitions of $1$ into two, and into three, such
elements. These would give $2t=1=3t$, a contradiction.
\ValidationDetails{%
\leanverified{BANANA/NonPrecompact/PrebananaHomogeneousLimitMassObstruction.lean}{conditional mass}
}%
\leanverifiedat{BANANA/NonPrecompact/BooleanBananaStoneNonamenability.lean}{100}{nonamenability of the actual limit automorphism group}
\ValidationDetails{%
\leanverifiedat{BANANA/NonPrecompact/BooleanInvariantMass.lean}{267}{finite additivity of the represented probability mass}
}%
\end{proof}

\subsection{\FolkBANANA}

\begin{proposition}\label{prop:folkman-no-eppa-amenability}
The class $\Kfolk$ does not have EPPA, and
$\Gfolk=\Aut(\Mfolk)$ is not amenable.
\end{proposition}
\begin{proof}
Use the zero right ring and a left ring with disjoint nonzero atoms $a,b$.
The isomorphism $\{0,a\}\to\{0,a\vee b\}$ cannot extend to a finite
ring automorphism, by the atom-count argument of
\Cref{prop:pre-banana-limits}.
\leanverifiedat{BANANA/NonPrecompact/FolkmanNativeEPPAObstruction.lean}{42}{failure of EPPA in the native two-sorted finite class}

For amenability, put $\Gfolk=\Aut(\Mfolk)$, fix a nonzero left element
$a$ of $\Mfolkset$ and let $U$ be its stabiliser. Then $U$ is an open subgroup of $\Gfolk$ and since amenability passes to open subgroups (see the proof of Lemma 13 in \cite{BPT}), it is enough to show that $U$ is not amenable to get that $\Gfolk$ is not amenable. 
\leanverifiedat{BANANA/NonPrecompact/OpenSubgroupAmenability.lean}{123}{amenability passes to open subgroups by coinduction}

The interval $[0,a]$ is a Boolean algebra on which $U$ acts by automorphisms.
Homogeneity, applied with the zero right ring, makes all $0<x<a$
conjugate under $U$. Partitions of $a$ into two and into three nonzero
pieces then yield the same contradiction $2t=1=3t$ on the Stone space of
$[0,a]$ as in \Cref{prop:pre-banana-limits}.
\ValidationRelative{\leanrepo/blob/\leancommit/BANANA/NonPrecompact/FolkmanFirstOrderLimitNonamenability.lean\#L34}%
  {nonamenability of the native paired automorphism group of the countable Fra\"iss\'e limit}
\ValidationDetails{%
\leanverifiedat{BANANA/NonPrecompact/FolkmanFirstOrderInterval.lean}{224}{actual left interval stabiliser transitivity}%
\leanverifiedat{BANANA/NonPrecompact/FolkmanFirstOrderInterval.lean}{301}{partition of any nonzero left element into any positive finite number of pieces}%
}%
\end{proof}

\subsection{BANANA: coherent EPPA and amenability}

\begin{theorem}\label{thm:eppa}
The BANANA class has coherent EPPA. 
A BANANA structure $\str A = (L_{A}, R_{A}, \beta_{\str A})$ has a coherent EPPA witness $\str B = (L_B, R_B, \beta_{\str B})$ with $\dim L_B = \dim R_B = \dim L_A + \dim R_A$. 
\end{theorem}
\leanverifiedat{BANANA/NonPrecompact/BananaNativeCoherentEppa.lean}{161}{coherent EPPA for arbitrary finite-dimensional pairings with the stated witness dimensions}
\ValidationDetails{%
  \leanverifiedat{BANANA/NonPrecompact/CoherentEppaWitness.lean}{374}{coherent EPPA in matrix coordinates}%
}%
\begin{proof}
Let $\str B=\widehat{\str A}$ be the perfect completion from \Cref{lem:completion}.
\leanverifiedat{BANANA/NonPrecompact/Completion.lean}%
  {80}{BananaMatrixStructure.completionEmbedding}
By \Cref{lem:perfect-homogeneous}, every partial automorphism extends to
$\str B$.
\leanverifiedat{BANANA/NonPrecompact/PerfectHomogeneity.lean}{157}{exists_perfectPairAutomorphism_extends}
We choose these extensions coherently.

For each isomorphism type of substructure of $\str A$, choose a representative
$\str U_0\le \str A$. Its perfect completion has side dimension at most $\dim L_A+\dim R_A$,
so it embeds in $\str B$. Start with any such embedding. Its restriction to $\str U_0$
and the given inclusion $\str U_0\le \str B$ are isomorphic finite substructures of
the perfect pairing $\str B$; by \Cref{lem:perfect-homogeneous}, an automorphism
of $\str B$ carries the former onto the latter pointwise. After composing, the
embedded completion therefore contains the given $\str U_0$ with the identity
inclusion.
\leanverifiedat{BANANA/NonPrecompact/AlignedCompletionEmbedding.lean}%
  {23}{BananaMatrixStructure.exists_alignedCompletionEmbedding}
An automorphism $(f,g)$ of $\str U_0$ acts on its completion by
\[
 f\oplus(g^{-1})^*\quad\text{on }L_{U_0}\oplus R_{U_0}^*,\qquad
 g\oplus(f^{-1})^*\quad\text{on }R_{U_0}\oplus L_{U_0}^*.
\]
This preserves the completed pairing and defines a group homomorphism.
\ValidationDetails{%
\leanverifiedat{BANANA/NonPrecompact/CompletionAutomorphism.lean}%
  {105}{BananaMatrixStructure.completionAutomorphism_left_completion}
\leanverifiedat{BANANA/NonPrecompact/CompletionAutomorphism.lean}%
  {121}{BananaMatrixStructure.completionAutomorphism_right_completion}
\leanverifiedat{BANANA/NonPrecompact/CompletionAutomorphismHom.lean}%
  {138}{BananaMatrixStructure.completionAutomorphism_comp_left}
\leanverifiedat{BANANA/NonPrecompact/CompletionAutomorphismHom.lean}%
  {161}{BananaMatrixStructure.completionAutomorphism_comp_right}
}%
A perfect subpair splits off orthogonally in a perfect ambient pairing:
on each sort its annihilator is a complementary subspace, and the two
complements pair perfectly.
\leanverifiedat{BANANA/NonPrecompact/BananaPerfectComplementPairing.lean}%
  {68}{left nondegenerate}
\leanverifiedat{BANANA/NonPrecompact/BananaPerfectComplementPairing.lean}%
  {104}{right nondegenerate}
Moreover every automorphism stabilising the
perfect subpair also stabilises its annihilator. We may therefore extend
the displayed action by the identity on the orthogonal complement. We obtain
a homomorphism
$\rho_{\str U_0}:\Aut(\str U_0)\to\Aut(\str B)$ extending every automorphism of $\str U_0$.
\leanverifiedat{BANANA/NonPrecompact/CoherentEppaPackage.lean}%
  {35}{BananaMatrixStructure.exists_alignedCoherentAmbientAction}
\ValidationDetails{%
\leanverifiedat{BANANA/NonPrecompact/EmbeddedCompletionAction.lean}%
  {79}{BananaMatrixStructure.embeddedCompletionAmbientAutomorphism_left}
\leanverifiedat{BANANA/NonPrecompact/EmbeddedCompletionAction.lean}%
  {97}{BananaMatrixStructure.embeddedCompletionAmbientAutomorphism_right}
\leanverifiedat{BANANA/NonPrecompact/EmbeddedCompletionAction.lean}%
  {146}{BananaMatrixStructure.embeddedCompletionAmbientAutomorphism_comp_left}
\leanverifiedat{BANANA/NonPrecompact/EmbeddedCompletionAction.lean}%
  {164}{BananaMatrixStructure.embeddedCompletionAmbientAutomorphism_comp_right}
}%

For each $\str U\le \str A$ isomorphic to $\str U_0$, choose $i_{\str U}:\str U_0\to \str U$ and an
automorphism $g_{\str U}\in\Aut(\str B)$ extending $i_{\str U}$. For $p:\str U\to \str V$ set
\[
 \Theta(p)=g_{\str V}\rho_{\str U_0}(i_{\str V}^{-1}p i_{\str U})g_{\str U}^{-1}.
\]
Indeed, on $\str U$ we have $g_{\str U}^{-1}=i_{\str U}^{-1}$, the middle factor restricts
to $i_{\str V}^{-1}pi_{\str U}$ on $\str U_0$, and $g_{\str V}=i_{\str V}$ on $\str U_0$; hence
$\Theta(p)|_{\str U}=p$.  For $q:\str V\to \str W$, cancellation of $g_{\str V}^{-1}g_{\str V}$ and
the homomorphism property of $\rho_{\str U_0}$ give
$\Theta(q)\Theta(p)=\Theta(qp)$, as required.
\leanverifiedat{BANANA/NonPrecompact/CoherentConjugation.lean}%
  {33}{conjugatedTransport_mul}
\end{proof}

\begin{proposition}\label{prop:amenable}
The group $G=\Aut(\str M)$ has a dense increasing union of finite subgroups.
In particular, $G$ is amenable.
\end{proposition}
\leanverifiedat{BANANA/NonPrecompact/BananaLimitTopologicalAmenability.lean}{196}{amenability of the pointwise automorphism group}
\begin{proof}
Exhaust $\mathbb Q$ by finite sets $Q_1\subseteq Q_2\subseteq\cdots$.
Let $K_n$ act on the perfect coordinate block $\Bban{Q_n}$ by
$T\in\mathrm{GL}(\F^{Q_n})$ on the left and $(T^{-1})^*$ on the right,
and act identically outside $Q_n$.  Extending $T$ by the identity on
$Q_{n+1}\setminus Q_n$ identifies $K_n$ with a subgroup of $K_{n+1}$;
each $K_n$ is finite.
\leanverifiedat{BANANA/NonPrecompact/FinitePerfectGL.lean}%
  {29}{finite GL stages}
\leanverifiedat{BANANA/NonPrecompact/FinitePerfectGL.lean}%
  {106}{injective block extension}
Given $g\in G$ and finitely many vectors, choose $Q_n$ containing their
supports and those of their images. The restriction of $g$ to the generated
substructure extends to an automorphism of $\Bban{Q_n}$ by
\Cref{lem:perfect-homogeneous}.
\ValidationDetails{%
\leanverifiedat{BANANA/NonPrecompact/PerfectHomogeneity.lean}{157}{exists_perfectPairAutomorphism_extends}
\leanverifiedat{BANANA/NonPrecompact/BananaLimitFiniteGLDensity.lean}%
  {25}{finite-set agreement}
\leanverifiedat{BANANA/NonPrecompact/BananaLimitDirectedGLUnion.lean}%
  {75}{common finite stage}
}%
\leanverifiedat{BANANA/NonPrecompact/BananaLimitLocallyFiniteDense.lean}%
  {24}{union agreement}
Thus $H=\bigcup_nK_n$ is a dense locally finite subgroup of $G$.
\leanverifiedat{BANANA/NonPrecompact/BananaLimitTopologicalAmenability.lean}{165}{increasing finite GL sequence with dense union}
Every locally finite group is amenable, and a topological group with a dense
amenable subgroup is amenable; see, for example,
\cite[Remark~4.3.6]{PestovDynamics}. Hence $G$ is amenable.
\end{proof}
\paragraph{Finite systems of partial automorphisms.}
Fix $n\ge1$. An $n$-\emph{system} is a tuple
$(\str A,p_1,\ldots,p_n)$ with $\str A\in\K$ and each $p_i$ being an isomorphism between
two substructures of $\str A$. An embedding of systems is an embedding of the
underlying structures which intertwines each $p_i$ with the corresponding
partial automorphism in the target. The class of $n$-systems has the
\emph{joint embedding property} (JEP) if any two systems embed into one.
It has the \emph{weak amalgamation property} (WAP) if every system $S$
has an extension $S'$ such that any two embeddings of $S'$ can be
amalgamated by maps which agree on the original copy of $S$; agreement on
all of $S'$ is not required. A subclass is \emph{cofinal} if every system
embeds into a member of it. Kechris--Rosendal prove that, for a Fra\"iss\'e
class, the diagonal conjugation action on the $n$th power of the limit
automorphism group has a comeagre orbit exactly when its class of
$n$-systems has JEP and WAP \cite[Theorem~6.2]{KR}.

\Needspace{6\baselineskip}
\begin{theorem}\label{thm:ample-generics}
The group $G$ has ample generics.
\end{theorem}
\begin{proof}

Fix $n$. Consider the subclass of $n$-systems in which the underlying
pairing is perfect and all $n$ named maps are total automorphisms.
We first show that this subclass has the amalgamation property and JEP, and
then that it is cofinal, which yields WAP and JEP for all $n$-systems. Suppose perfect total systems
$\str B$ and $\str C$ share a perfect total subsystem $\str A$. Because the restriction
of the pairing to $\str A$ is perfect, we obtain orthogonal
decompositions $
 \str B=\str A\perp \str B'$ and $\str C=\str A\perp \str C'$ by setting
\[
 L_{B'}=\{x\in L_B:\beta_{\str B}(x,R_A)=0\},\qquad
 R_{B'}=\{y\in R_B:\beta_{\str B}(L_A,y)=0\},
\]

and likewise for $\str C'$.
\leanverifiedat{BANANA/NonPrecompact/BananaPerfectOrthogonalEquiv.lean}%
  {62}{orthogonal splitting}
\leanverifiedat{BANANA/NonPrecompact/BananaPerfectOrthogonalEquiv.lean}%
  {72}{dual orthogonal splitting}
Every named automorphism preserves $\str A$ and hence these annihilators.
\ValidationDetails{%
\leanverifiedat{BANANA/NonPrecompact/PerfectTotalInvariantComplements.lean}%
  {63}{left annihilator invariant}
\leanverifiedat{BANANA/NonPrecompact/PerfectTotalInvariantComplements.lean}%
  {79}{right annihilator invariant}
\leanverifiedat{BANANA/NonPrecompact/PerfectTotalComplementEquiv.lean}%
  {25}{left complement action}
\leanverifiedat{BANANA/NonPrecompact/PerfectTotalComplementEquiv.lean}%
  {68}{right complement action}
}%
On $ \str D=\str A\perp \str B'\perp \str C'$
combine the corresponding automorphisms on the three summands. Their
restrictions agree on $\str A$, so this simultaneously amalgamates all $n$
named automorphisms. Orthogonal sums over the zero structure give joint
embedding in the same subclass.
\leanverifiedat{BANANA/NonPrecompact/BananaPerfectTotalSplitSystems.lean}%
  {42}{split-system amalgamation}
\leanpartial{BANANA/NonPrecompact/BananaPerfectTotalDiagramAmalgam.lean}%
  {perfect-AP}
\ValidationDetails{%
  \leanverifiedat{BANANA/NonPrecompact/BananaNativePerfectSystemAmalgam.lean}{115}{one native amalgam for all named automorphisms of perfect systems}%
}%

As stated above, this subclass is also cofinal. Given an $n$-system
$(\str A,p_1,\ldots,p_n)$, apply \Cref{lem:completion} to $\str A$ to
embed it in a perfect pairing,
\leanverifiedat{BANANA/NonPrecompact/Completion.lean}%
  {80}{BananaMatrixStructure.completionEmbedding}
and then use \Cref{lem:perfect-homogeneous} to extend each $p_i$ to the same perfect completion.
\leanverifiedat{BANANA/NonPrecompact/PerfectHomogeneity.lean}{157}{exists_perfectPairAutomorphism_extends}
\leanverifiedat{BANANA/NonPrecompact/BananaFiniteSystemPerfectTotalCofinality.lean}%
  {40}{perfect-total cofinality}
\ValidationDetails{%
  \leanverifiedat{BANANA/NonPrecompact/BananaNativePerfectSystems.lean}{69}{a common native perfect completion for all named maps}%
}%

Since this subclass is cofinal and has the amalgamation property and JEP, the
full class of $n$-systems has WAP and JEP (enlarging a system and any two
targets into the perfect-total subclass allows us to amalgamate or jointly
embed them there).
\leanverifiedat{BANANA/NonPrecompact/BananaFiniteSystemsJepWap.lean}%
  {167}{system JEP/WAP}
\ValidationDetails{%
\leanverifiedat{BANANA/NonPrecompact/BananaLimitDenseConjugacy.lean}{187}{a dense diagonal conjugacy orbit in each finite power}%
\leanverifiedat{BANANA/NonPrecompact/BananaLimitDenseConjugacy.lean}{196}{the set of points with dense diagonal conjugacy orbit is comeagre}%
}%
The Kechris--Rosendal criterion therefore gives a comeagre diagonal
conjugacy orbit in $G^n$. Since $n$ was arbitrary, $G$ has ample generics.
\leanverifiedat{BANANA/NonPrecompact/BananaLimitGenericBackAndForth.lean}{285}{a single comeagre diagonal conjugacy orbit in every finite power}
\ValidationDetails{%
\leanverifiedat{BANANA/NonPrecompact/BananaLimitGenericExtension.lean}{273}{the perfect-system coverage and extension conditions define a residual set of tuples}%
\leanverifiedat{BANANA/NonPrecompact/BananaLimitGenericBackAndForth.lean}{220}{conjugacy from finite perfect-system extensions by back-and-forth}%
}%
\end{proof}

\subsection{Definable versus unrestricted Ramsey colourings}
\label{sec:definable-ramsey}

Following Hrushovski--Krupi\'nski--Pillay
\cite[Definition~2.9]{HKPFirstOrder}, a complete theory $T$ is
\emph{first-order extremely amenable} if every complete type over
$\emptyset$, in any tuple of variables, has a global complete extension
invariant under $\Aut(\mathfrak C)$, where
$\mathfrak C\models T$ is a sufficiently saturated and homogeneous model.
This is a property of the theory and is distinct from extreme amenability of
the topological group $\Aut(\str M)$.

Following \cite{KLMRamseyDynamics}, $T$ has the \emph{definable elementary
embedding Ramsey property} (DEERP) if, for finite tuples $\bar a\subseteq\bar b$
in $\mathfrak C$, every finite colouring of the realisations of
$\operatorname{tp}(\bar a/\emptyset)$ with fibres definable with parameters
in $\mathfrak C$ admits a realisation $\bar b'$ of
$\operatorname{tp}(\bar b/\emptyset)$ such that all such realisations of
$\operatorname{tp}(\bar a/\emptyset)$ with entries in $\bar b'$ have the same colour.

\begin{corollary}[Definable versus unrestricted Ramsey colourings]
\label{cor:deerp-separation}
Let $\str M$ be the BANANA limit and put $T=\operatorname{Th}(\str M)$.  Then $T$
has the definable elementary embedding Ramsey property (DEERP) of
Krupi\'nski--Lee--Moconja.  Nevertheless, every finite
$\str A=(L_A,R_A,\beta_{\str A})$ with $L_A\ne0$ and $R_A\ne0$ has infinite ordinary
copy Ramsey degree in $\operatorname{Age}(\str M)$.
\end{corollary}
\begin{proof}
The BANANA class has canonical amalgamation over $\emptyset$ in the sense
of Hrushovski--Krupi\'nski--Pillay; here the substructure generated by
$\emptyset$ consists of the two named zero elements.  For two structures
$\str A$ and $\str B$, take the direct sums $L_A\oplus L_B$ and $R_A\oplus R_B$,
retain the two old pairings, and put both cross-pairings equal to zero.
\ValidationDetails{%
\leanpartialat{BANANA/NonPrecompact/BananaDirectSum.lean}%
  {181}{BananaMatrixStructure.directSumInlEmbedding}
\leanpartialat{BANANA/NonPrecompact/BananaDirectSum.lean}%
  {199}{BananaMatrixStructure.directSum_cross_pairing_left_right}
\leanverifiedat{BANANA/NonPrecompact/BananaDirectSum.lean}%
  {367}{BananaMatrixStructure.directSumEmbedding}
\leanverifiedat{BANANA/NonPrecompact/BananaDirectSum.lean}%
  {424}{BananaMatrixStructure.directSumEmbedding_comp_left}
\leanverifiedat{BANANA/NonPrecompact/BananaDirectSum.lean}%
  {447}{BananaMatrixStructure.directSumEmbedding_comp_right}
}%
\leanverified{BANANA/NonPrecompact/BananaDirectSumCanonical.lean}{zero-base AP}
This construction is functorial for embeddings: an embedding on either
factor extends by the identity on the other factor (and similarly for
embeddings on both factors). Since the limit is $\omega$-categorical by
\Cref{prop:fraisse}, \cite[Proposition~2.17]{HKPFirstOrder} implies that
$T$ is first-order extremely amenable. By
\cite[Theorem~3.16(ii)]{KLMRamseyDynamics}, this is equivalent to DEERP.
\ValidationRelative{\leanrepo/blob/\leancommit/BANANA/NonPrecompact/BananaLimitDefinableColouring.lean\#L190}%
  {definable elementary embedding Ramsey property for the complete tagged relational BANANA theory}
\ValidationDetails{%
  \ValidationRelative{\leanrepo/blob/\leancommit/BANANA/NonPrecompact/BananaLimitInvariantTypes.lean\#L34}%
  {every complete type over the empty set of the tagged BANANA theory has an invariant complete extension}%
  \leanverifiedat{BANANA/NonPrecompact/FirstOrderCanonicalAmalgamation.lean}{163}{invariant extensions from canonical age amalgamation for countable ultrahomogeneous omega-categorical structures}%
  \leanverifiedat{BANANA/NonPrecompact/FirstOrderInvariantEnumerationRamseySaturated.lean}{136}{definable tuple Ramsey is equivalent to an invariant enumeration extension under finite-parameter saturation and strong finite homogeneity}%
  \leanverifiedat{BANANA/NonPrecompact/FirstOrderEnumerationEllis.lean}{416}{definable tuple Ramsey is equivalent to the Ellis-map condition under finite-parameter saturation and strong finite homogeneity}%
  \leanverifiedat{BANANA/NonPrecompact/BananaLimitEnumerationEllis.lean}{43}{an Ellis map with invariant values on realised enumeration types of the tagged BANANA limit}%
  \leanverifiedat{BANANA/NonPrecompact/BananaSortedDefinableTupleRamsey.lean}{169}{definable tuple Ramsey in the two-sorted language for every model decoded from the tagged BANANA theory}%
}%
The second assertion is \Cref{thm:overview}(iii).
\ValidationDetails{%
\ValidationRelative{\leanrepo/blob/\leancommit/BANANA/NonPrecompact/BananaLimitParameterColouring.lean\#L106}%
  {monochromatic finite targets for colourings invariant under finite parameter stabilisers in the concrete limit}%
}%
\end{proof}

In particular, $\Aut(\str M)$ itself is not extremely amenable; in fact
\Cref{thm:dynamics} proves the stronger statement that it contains no
co-precompact extremely amenable subgroup.

Thus every definable finite colouring in the sense of
Krupi\'nski--Lee--Moconja can be made monochromatic on a copy of any
prescribed finite configuration, whereas unrestricted finite colourings of
copies of a mixed BANANA structure cannot be reduced to any fixed finite Ramsey degree.
For the two rigid four-element line pairs $\Aban{0},\Aban{1}$ used in
\Cref{sec:persistent-colouring}, copy and embedding colourings coincide, so
the same contrast already holds for unrestricted embedding colourings.  In
particular, the persistent integer support/intersection colourings used below
cannot witness a failure of DEERP. The point is not merely that first-order
extreme amenability and extreme amenability of the topological automorphism
group differ---examples of that distinction were already known---but that
BANANA places the distinction in a locally finite,
$\omega$-categorical and amenable Fra\"iss\'e limit while the ordinary
Ramsey degrees remain infinite and the universal minimal flow is nonmetrisable.

\section{Non-existence of precompact Ramsey expansions}
\label{sec:no-precompact}
We show that none of the three classes admits a precompact Ramsey expansion.
For each, we colour copies of a fixed small structure so that every
copy of a suitable finite target contains arbitrarily many colours.
\subsection{\BoolBANANA}
\label{sec:pre-banana-colouring}

We use the following elementary lemma.

\begin{lemma}\label{lem:pre-banana-subset-sums}
If $q\ge2$ is a power of $2$, the subset sums of any $q-1$ odd integers
cover $\mathbb Z/q\mathbb Z$.
\end{lemma}
\leanverifiedat{BANANA/NonPrecompact/IndexedSubsetSums.lean}{224}{odd_indexedSubsetSums_two_pow_eq_univ}
\begin{proof}
Starting with $S=\{0\}$, we add the odd integers one after the other by repeatedly replacing $S$ with
$S\cup(S+a)$. 
Unless $S$ is the whole group $\mathbb{Z}/q \mathbb{Z}
$, this increases the size of $S$, as otherwise $S$ would be invariant under translation by $a$, which generates
$\mathbb Z/q\mathbb Z$. 
Thus after $q-1$ steps all residues occur.
\end{proof}

\begin{theorem}\label{thm:pre-banana-colouring}
Let $\Abool\in\Kbool$ have two atoms, both marked. Let $q \geq 2$ be a power of $2$. 
Let $\Bbool{q}\in\Kbool$ have $q$ atoms, all marked.
Every $\str C\in\Kbool$ admits a $q/2$-colouring of its $\Abool$-copies
such that every $\Bbool{q}$-copy contains all $q/2$ colours.
Consequently $\Abool$ has infinite copy Ramsey degree and $\Kbool$
has no precompact Ramsey expansion.
\end{theorem}
\ValidationRelative{\leanrepo/blob/\leancommit/BANANA/NonPrecompact/BooleanBananaCopyTransport.lean\#L274}%
  {persistent colouring on native Boolean copy ranges of the canonical source}
\begin{proof}
We pick an atom $c_*$ of $\str C$. 
A copy of $\Abool$ in $\str C$ is an unordered partition of $\At(\str C)$ into two blocks with odd numbers
of marked atoms. 
We colour the copy by the residue of the number of marked atoms in the block containing $c_*$ modulo $q$. 
This is well-defined on copies, despite the automorphism interchanging the two atoms of $\Abool$.

Fix any copy $\str B$ of $\Bbool{q}$ in $\str C$, with blocks $V_1,\ldots,V_q$ labelled such that
$c_*\in V_1$. 
Let $a_i$ be the number of marked atoms in $V_i$, so
each $a_i$ is odd. For $S\subseteq\{2,\ldots,q\}$ the block
$V_1\cup\bigcup_{i\in S}V_i$
contains $a_1+\sum_{i\in S}a_i$ marked atoms. 
By \Cref{lem:pre-banana-subset-sums}, these values cover all residues modulo $q$. 
An odd residue forces $|S|$ to be even. 
The chosen block and its complement then each contain an odd number of the $V_i$. 
In particular, they are nonempty
and contain an odd number of marked atoms.
\leanverifiedat{BANANA/NonPrecompact/ColouringWrappers.lean}{32}{prebanana_odd_residue_copy_witness}
Thus they form a copy of $\Abool$ of the required
colour inside $\str B$. 
Since $q/2$ can be arbitrarily large, the copy Ramsey degree of $\Abool$ is infinite.
By \Cref{lem:expansion-degree-bound}, $\Kbool$ does not have a precompact Ramsey expansion.
\ValidationDetails{%
\leanverifiedat{BANANA/NonPrecompact/PrebananaPersistence.lean}%
  {124}{exists_prebananaPersistentColouring}
\leanverifiedat{BANANA/NonPrecompact/PrebananaCopyDegree.lean}%
  {125}{prebananaTwoMarkedCopyRamseyDegree_infinite}
}%
\leanverifiedat{BANANA/NonPrecompact/BananaGeneratedRamseyObstructions.lean}{21}{no precompact Ramsey expansion of the marked Boolean age}
\ValidationDetails{%
  \ValidationRelative{\leanrepo/blob/\leancommit/BANANA/NonPrecompact/NativePointedRamseyExpansion.lean\#L119}%
  {native Boolean--BANANA obstruction via the generated pointed cover}%
}%
\end{proof}

Note that the atom $c_*$ is not named in the structures and does not have to be preserved
by embeddings. It is only used to select one of the two blocks to define the colouring.

\subsection{\FolkBANANA}
\label{sec:folkman-colouring}

Let $\Afolk\in\Kfolk$ be a \FolkBANANA{} structure with one atom on each side with pairing one. Let $q \geq 2$ be a power of~$2$.
Let $\Bfolk{q} = (L, R, \beta)$ be the \FolkBANANA{} structure with $\At(L) =  \{l_1,\ldots,l_q\}$, $\At(R) = \{r_1, \dots, r_q\}$ and $\beta(l_i,r_j)=1$ exactly
when $i=j$. 
Let $\str C = (L_C, R_C, \beta_{\str C}) \in \Kfolk$ and define
\begin{equation}\label{eq:folkman-integer-count}
 N_{\str C}(x,y)=\bigl|\{(u,v)\in\At(L_C)\times\At(R_C):
 u\le x,\ v\le y,\ \beta_{\str C}(u,v)=1\}\bigr|.
\end{equation}
Additivity gives $N_{\str C}(x,y)\equiv\beta_{\str C}(x,y)\pmod2$.

\begin{theorem}\label{thm:folkman-infinite-degree}
Every $\str C\in\Kfolk$ admits a $q/2$-colouring of its
$\Afolk$-copies such that every $\Bfolk{q}$-copy contains all
$q/2$ colours. Thus $\Afolk$ has infinite copy Ramsey degree.
In particular, the class $\Kfolk$ has no precompact Ramsey expansion.
\end{theorem}
\ValidationRelative{\leanrepo/blob/\leancommit/BANANA/NonPrecompact/FolkmanRamseyExpansionObstruction.lean\#L114}%
  {persistent colouring on native copy ranges of the canonical source}
\ValidationDetails{%
\leanverifiedat{BANANA/NonPrecompact/FolkmanRamseyExpansionObstruction.lean}{130}{infinite native copy Ramsey degree}
\leanpartialat{BANANA/NonPrecompact/FolkmanPersistence.lean}{240}{atom-code colouring}
\leanpartialat{BANANA/NonPrecompact/FolkmanPersistence.lean}{379}{atom-code copy degree}
}%
\begin{proof}
A copy of $\Afolk$ in $\str C$ is determined by nonzero $x \in L_C$ and $y \in R_C$ such that $\beta_{\str C}(x,y)=1$.
We colour the copy by the residue of $N_{\str C}(x,y)$ modulo $q$. 
Thus the possible colours are the odd residues modulo $q$.
In a $\Bfolk{q}$-copy $\str B$ in $\str C$, write $x_i\in L_C$ and $y_i\in R_C$ for the atom images and put
$y=y_1\vee\cdots\vee y_q$. Since the $y_j$ are pairwise disjoint, each $a_i=N_{\str C}(x_i,y)$ is odd because
$\beta_{\str C}(x_i,y)=\sum_j\beta_{\str C}(x_i,y_j)=1$.
Since the $x_i$ are pairwise disjoint,
\[
 N_{\str C}\Bigl(\bigvee_{i\in S}x_i,y\Bigr)=\sum_{i\in S}a_i.
\]
Apply \Cref{lem:pre-banana-subset-sums} to the first $q-1$ of these
integers. Each odd residue $r$ modulo $q$ is attained by a nonempty set $S_r$, and for $x=\bigvee_{i\in S_r}x_i$ the integer $N_{\str C}(x,y)$ is odd. 
Hence the resulting left element $x$
is nonzero and $\beta_{\str C}(x,y) =1$, so $x$ and $y$ define a copy of $\Afolk$ in $\str B$ with colour $r$.
\leanverifiedat{BANANA/NonPrecompact/ColouringWrappers.lean}{95}{folkman_odd_residue_copy_witness}
Only joins of atom images belonging to the target copy are used.
Atoms outside the image of the non-unital embedding play no role.

By Lemma~\ref{lem:expansion-degree-bound}, $\Kfolk$ has no precompact Ramsey expansion.
\ValidationRelative{\leanrepo/blob/\leancommit/BANANA/NonPrecompact/NativePointedRamseyExpansion.lean\#L125}%
  {native Folkman--BANANA obstruction via the generated pointed cover}
\ValidationDetails{%
  \leanverifiedat{BANANA/NonPrecompact/BananaGeneratedRamseyObstructions.lean}{33}{no precompact Ramsey expansion of the encoded Folkman product class}%
}%
\end{proof}

\subsection{BANANA}
\label{sec:persistent-colouring}

For $b\in\{0,1\}$, let $\Aban{b} = (L_{A_b},R_{A_b},\beta_{\Aban{b}})$ with $L_{A_b} = \langle x\rangle$, $R_{A_b} = \langle  y \rangle$ being one-dimensional, and $\beta_{\Aban{b}}(x,y) = b$.
Both $\Aban{0}$ and $\Aban{1}$ are rigid.
Fix $\str C\in\K$ and let $q\ge2$ be a power of $2$.
By \Cref{lem:completion}, choose a coordinate realisation $\str C\hookrightarrow
\Bban{N}$ and identify coordinate vectors with subsets of $[N]$.
A copy of $\Aban{b}$ in $\str C$ is determined by the unique nonzero left and right vectors $x$ and $y$, and we colour the copy with the colour
\begin{equation}\label{eq:intersection-colouring}
|x\cap y|\pmod q.
\end{equation}
Clearly, only the $q/2$ residues of parity $b$ can occur.
\leanverifiedat{BANANA/NonPrecompact/CompletionPersistence.lean}%
  {123}{completionIntersectionCount_map_parity}
\leanverifiedat{BANANA/NonPrecompact/CompletionPersistence.lean}%
  {74}{residueParityPalette_card}

Put $D=q-1$. 
We work in the quotient ring
$\mathcal R_q=\F[z]/(z^q-1)$, so exponents of~$z$ are read modulo $q$.
\ValidationRelative{\leanrepo/blob/\leancommit/BANANA/NonPrecompact/ResidueAlgebra.lean\#L22}%
  {ResidueAlgebra}
Every element of $\mathcal R_q$ has a unique representative of the form
$\sum_{\rho=0}^{q-1}\alpha_\rho z^\rho$ with $\alpha_\rho \in \F$.
For nonnegative integers $(M_\omega)$, the coefficient of $z^\rho$ in
$\sum_\omega z^{M_\omega}$ is the parity of the number of indices
with $M_\omega\equiv\rho\pmod q$.
Setting $t=1+z$ gives
\begin{equation}\label{eq:sum-lemma-ring}
 t^q=0,\qquad t^D=1+z+\cdots+z^{q-1}.
\end{equation}
Indeed $(1+z)^q=1+z^q$ in characteristic~$2$, while
$(1+z)^{q-1}=\prod_{i=0}^{\log_2q-1}(1+z^{2^i})$.
An embedding $\Bban{D}\hookrightarrow \Bban{N}$ is given by matrices
$P,Q\in\F^{N\times D}$ with
\begin{equation}\label{eq:PQ-identity}
 P^{\mathsf T}Q=I_D.
\end{equation}
Writing their rows as $a_e^{\mathsf T},b_e^{\mathsf T}$, the size of the intersection of the images of $u,v\in\F^D$ in $\Bban{N}$ is given by the integer
\[
 F(u,v)=\sum_{e=1}^N[a_e\cdot u]\,[b_e\cdot v],\qquad
 \sum_{e=1}^N a_e b_e^{\mathsf T}=I_D.
\]
Its parity is $u\cdot v$. The next statement allows affine coordinate
functions as well, so that later we can ensure that both vectors are nonzero.

\begin{proposition}
\label{prop:residue-surjectivity}\label{prop:sharp-parity-residues}
Let $q\ge2$ be a power of $2$ and $D=q-1$. For a finite set $E$, choose
$a_e,b_e\in\F^D$, $c_e,d_e\in\F$ with
$\sum_{e\in E}a_e b_e^{\mathsf T}=I_D$. Set
$f_e(u)=[a_e\cdot u+c_e]$, $g_e(v)=[b_e\cdot v+d_e]$ and $
 F(u,v)=\sum_{e\in E}f_e(u)g_e(v)$.
Then $F$ attains every residue modulo $q$ an odd number of times on
$\F^D\times\F^D$.
\end{proposition}
\leanverifiedat{BANANA/NonPrecompact/BilinearResidue.lean}%
  {533}{bilinearResidueFiber_card_odd}
\begin{proof}
In $\mathcal R_q$, each bit $f_e(u)g_e(v)$ satisfies
$z^{f_e(u)g_e(v)}=1+t f_e(u)g_e(v)$. Thus, for each $(u,v)$,
\[
 z^{F(u,v)}=\prod_{e\in E}z^{f_e(u)g_e(v)}
 =\prod_{e\in E}\bigl(1+t f_e(u)g_e(v)\bigr).
\]
Expanding this product and then summing over $u,v$ gives
\begin{equation}\label{eq:bilinear-sum-expansion}
 \sum_{u,v}z^{F(u,v)}=\sum_{S\subseteq E}t^{|S|}N(S),\qquad
 N(S)=\sum_{u,v}\prod_{e\in S}f_e(u)g_e(v).
\end{equation}
Here $N(S)$ is taken modulo two and counts solutions of the system of equations
$a_e\cdot u+c_e=1, \, b_e\cdot v+d_e=1$ for $e \in S$.
If $|S|<D$, there are fewer than $D$ affine equations in the
$D$ coordinates of $u$, and independently fewer than $D$ in those of
$v$. 
Each system has either no solutions or an even number of solutions:
if it is consistent, its solution set is a translate of the kernel of a
linear map $\F^D\to\F^S$; since $|S|<D$, that kernel has positive
dimension and therefore even cardinality. Thus $N(S)=0$.
If $|S|>D$, then $t^{|S|}=0$ by
\eqref{eq:sum-lemma-ring}.

For $|S|=D$, let $A_S,B_S$ have columns $(a_e)_{e\in S}$ and
$(b_e)_{e\in S}$ in the same order. Each system has an odd number of
solutions precisely when its matrix is invertible, in which case its
solution is unique. Thus $N(S)=\det(A_S)\det(B_S)$ over $\F$.
Let $A,B$ consist of all these columns. Cauchy--Binet and
$AB^{\mathsf T}=I_D$ now give
\[
 \sum_{u,v}z^{F(u,v)}
 =t^D\sum_{|S|=D}\det(A_S)\det(B_S)
 =t^D\det(AB^{\mathsf T})
 =1+z+\cdots+z^{q-1}.
\]
Comparing coefficients proves the assertion.
\end{proof}

\begin{theorem}
\label{thm:four-element-infinite-degree}\label{thm:infinite-degree}
Let $q \geq 2$ be a power of $2$. 
Each $\str C\in\K$ has a $q/2$-colouring of
its $\Aban{b}$-copies such that every $\Bban{q}$-copy contains all $q/2$ colours,
for each $b\in\{0,1\}$. For $\Aban{1}$, it already suffices to consider $\Bban{q-1}$.
Consequently, both $\Aban{0}$ and $\Aban{1}$ have infinite copy Ramsey degree.
\end{theorem}
\leaninterfaceat{BANANA/NonPrecompact/CopyRamseyDegree.lean}%
  {162}{exists_completionResiduePersistentColouring}
\leaninterfaceat{BANANA/NonPrecompact/CopyRamseyDegree.lean}%
  {236}{exists_pairingOne_completionResiduePersistentColouring}
\leaninterfaceat{BANANA/NonPrecompact/CopyRamseyDegree.lean}%
  {328}{linePairCopyRamseyDegree_infinite}
\begin{proof}
Use the coordinate realisation and colouring in
\eqref{eq:intersection-colouring}. For $\Aban{1}$, take any copy of $\Bban{q-1}$ in $\str C$.
Its matrices satisfy \eqref{eq:PQ-identity}, so
\Cref{prop:residue-surjectivity}, with $c_e=d_e=0$, gives every residue.
An odd residue is realised by two vectors $x,y \neq 0$ with $\beta_{\Bban{q-1}}(x,y) = 1$, i.e., by a copy of $\Aban{1}$.
\leanverifiedat{BANANA/NonPrecompact/CompletionPersistence.lean}%
  {227}{exists_pairingOne_completionIntersectionColour}

For either value of $b$, take a $\Bban{q}$-copy and restrict to
$x=(u,1)$, $y=(v,1)$, with $u,v\in\F^{q-1}$.
Both vectors are nonzero. The ambient coordinates become affine functions
$f_e(u),g_e(v)$, and pairing preservation gives
$
 \sum_e f_e(u)g_e(v)\equiv u\cdot v+1\pmod2$.
Comparing mixed coefficients gives
$\sum_e a_e b_e^{\mathsf T}=I_{q-1}$. Equivalently, this coefficient
identity is the upper-left $(q-1)\times(q-1)$ block of
$P^{\mathsf T}Q=I_q$ for the chosen $\Bban{q}$-copy. The affine residue count therefore
gives every residue on this set of pairs. A residue of parity $b$
corresponds to an $\Aban{b}$-copy.
\leanverifiedat{BANANA/NonPrecompact/CompletionPersistence.lean}%
  {205}{exists_completionIntersectionPalette}
The same ambient colouring works for every target copy, and $q/2$ is
unbounded.
\end{proof}

For this particular colouring, dimension $q-1$ is necessary for the odd
source: in a standard coordinate block of smaller dimension the colour
$q-1$ cannot occur. 

\subsection{All Ramsey degrees}
\label{sec:degree-classification}

\begin{lemma}\label{lem:degree-propagation}
Let $\mathcal C$ be an amalgamation class and let $\str P,\str A\in\mathcal C$
with $\str P$ embeddable in $\str A$.
If $s=|\copies{\str A}{\str P}|$ and $t^{\rm cp}_{\mathcal C}(\str A)=d<\infty$, then
$t^{\rm cp}_{\mathcal C}(\str P)\le sd$.
\end{lemma}
\leanverifiedsuccessorat{SuccessorTree/CopyRamsey.lean}%
  {257}{copyRamseyDegreeLE_of_embedding}
\begin{proof}
Fix $\str B\in\mathcal C$ and a number of colours $r$.
Successively amalgamate a copy of $\str A$ over each $\str P$-copy of $\str B$, carrying
forward the copies attached at earlier stages, and obtain $\str B^+$ in which
every such copy extends to an $\str A$-copy.
Choose a degree-$d$ witness for source $\str A$, target $\str B^+$ and $2^r$ colours.
Given an $r$-colouring of $\str P$-copies, colour each $\str A$-copy by the set of
colours occurring inside it, a set of size at most $s$.
Some $\str B^+$-copy sees at most $d$ such sets. Every $\str P$-copy in its
distinguished $\str B$-copy extends to one of these $\str A$-copies, so at most
$sd$ colours occur in that $\str B$-copy. If $\str B$ has no $\str P$-copy, the
required bound is vacuous.
\end{proof}

\begin{theorem}\label{thm:banana-degree-classification}
For every $\str A=(L_A,R_A,\beta_{\str A})\in\K$,
\[
 t^{\rm cp}_{\K}(\str A)=
 \begin{cases}
  1,&L_A=0\text{ or }R_A=0,\\
  \infty,&L_A\ne0\text{ and }R_A\ne0,
 \end{cases}
\]
and consequently $t_{\K}(\str A)=|\Aut(\str A)|$ if $L_A=0$ or $R_A=0$, and $t_{\K}(\str A)=\infty$ otherwise.
\end{theorem}
\ValidationRelative{\leanrepo/blob/\leancommit/BANANA/NonPrecompact/RamseyDegreeClassifications.lean\#L15}%
  {unconditional BANANA embedding-degree classification in matrix coordinates}
\ValidationDetails{%
\leanverifiedat{BANANA/NonPrecompact/BananaOneSidedCopyDegree.lean}%
  {157}{copy degrees}
}%
\begin{proof}
Fix a number of colours $r\ge1$.
If $R_A=0$ and $a=\dim L_A>0$, we use \Cref{lem:completion} to choose
$D\ge a$ and an embedding of a given target $\str B$ into $\Bban{D}$. We use the
following special case of the finite vector-space
Ramsey theorem \cite{GLR}: for $a,D,r$ there is $N$ such that every
$r$-colouring of the $a$-dimensional subspaces of $\F^N$ has a
$D$-dimensional subspace all of whose $a$-dimensional subspaces have one
colour. Choose such $N$ and call the monochromatic $D$-space $W$.
\leanverified{BANANA/NonPrecompact/LinearBinaryRepresentability.lean}{binary GLR}
Use the ambient pairing $(V,V^*,\mathrm{evaluation})$. Its $\str A$-copies
are precisely those left subspaces with right side zero.
Extend the dual coordinate functionals of a basis of $W$ to $V$.
Their span pairs perfectly with $W$, giving a $\Bban{D}$-copy, and hence a
$\str B$-copy with all its $\str A$-copies monochromatic.
\ValidationDetails{%
\leanverifiedat{BANANA/NonPrecompact/PerfectLeftRange.lean}%
  {37}{BananaMatrixStructure.exists_perfectBlockEmbedding_leftRange_eq}
\leaninterfaceat{BANANA/NonPrecompact/LinearBinaryRepresentability.lean}%
  {108}{leftOneSidedCopyRamseyOne_via_successors}
}%
If both sides of $\str A$ are zero there is exactly one copy in each structure;
the case $L_A=0$ is symmetric.
\ValidationDetails{%
\leanverifiedat{BANANA/NonPrecompact/PerfectLeftRange.lean}%
  {117}{BananaMatrixStructure.exists_perfectBlockEmbedding_rightRange_eq}
\leaninterfaceat{BANANA/NonPrecompact/LinearBinaryRepresentability.lean}%
  {114}{rightOneSidedCopyRamseyOne_via_successors}
}%

If both sides are nonzero, choose nonzero $x\in L_A,y\in R_A$.
Their lines form an $\Aban{b}$-copy for $b=\beta_{\str A}(x,y)$.
\ValidationDetails{%
\leanverifiedat{BANANA/NonPrecompact/TwoSidedLinePair.lean}%
  {49}{line-pair exists}
\leanverifiedat{BANANA/NonPrecompact/LinePairSourceEmbedding.lean}%
  {71}{line-pair copy}
\leaninterfaceat{BANANA/NonPrecompact/CopyRamseyDegree.lean}%
  {328}{linePairCopyRamseyDegree_infinite}
}%
The claim for $t^{\rm cp}_{\K}(\str A)$ then follows from \Cref{thm:four-element-infinite-degree,lem:degree-propagation}, and $t_{\K}(\str A)=|\Aut(\str A)|\,t^{\rm cp}_{\K}(\str A)$ gives the corresponding values of $t_{\K}(\str A)$.
\leanverifiedat{BANANA/NonPrecompact/BananaCopyDegree.lean}%
  {170}{mixed copies}
\leanverifiedsuccessorat{SuccessorTree/CopyRamsey.lean}%
  {257}{copyRamseyDegreeLE_of_embedding}
\end{proof}

Thus the structures $\Aban{b}$ with $b \in \{0,1\}$ are the smallest BANANA structures with infinite Ramsey degree.

\begin{corollary}\label{cor:no-precompact}\label{cor:questions}
BANANA has no precompact Ramsey expansion. Consequently
Conjecture~\ref{conj:eppa-ramsey} is false even with coherent EPPA and an
$\omega$-categorical limit.
\end{corollary}
\begin{proof}
Combine \Cref{lem:expansion-degree-bound,thm:infinite-degree,prop:fraisse,thm:eppa}.
\ValidationRelative{\leanrepo/blob/\leancommit/BANANA/NonPrecompact/NativePointedRamseyExpansion.lean\#L132}%
  {native BANANA obstruction for hereditary expansions via the generated pointed cover}
\ValidationDetails{%
  \leanverifiedat{BANANA/NonPrecompact/BananaGeneratedRamseyObstructions.lean}{44}{no precompact Ramsey expansion of the tagged functional BANANA age}%
}%
\end{proof}

\subsection{Dynamical consequences}

\begin{lemma}\label{lem:EA-degree}
If the automorphism group of a countable locally finite homogeneous
structure $\str N$ has a co-precompact extremely amenable subgroup, every
member of $\Age(\str N)$ has finite copy Ramsey degree.
\end{lemma}
\begin{proof}

The finite case is immediate. Taking the closure of the subgroup
preserves extreme amenability and co-precompactness. The proof of
\cite[Corollary~3.3]{BMT} gives metrisability of $M(\Aut(\str N))$.
By \cite[Theorem~8.7]{ZuckerDynamics16}, the embedding Ramsey degrees are
finite, and these bound the copy Ramsey degrees.

\end{proof}
\ValidationRelative{\leanrepo/blob/\leancommit/BANANA/NonPrecompact/FirstOrderExtremelyAmenableDegree.lean\#L175}%
  {direct proof of finite copy Ramsey degree from co-precompact extreme amenability in the first-order interface}
\ValidationDetails{%
\leanpartialat{BANANA/NonPrecompact/FixedPointColouring.lean}{144}{finite palette from a fixed point and finitely many orbit representatives}
  \leanverifiedat{BANANA/NonPrecompact/CoprecompactFixedPointMetrizability.lean}{118}{metrisability of a minimal compact flow from a co-precompact extremely amenable subgroup}%
  \leanverifiedat{BANANA/NonPrecompact/FirstOrderMetrizableFlowRamseyDegree.lean}{76}{finite copy Ramsey degrees from a metrisable universal minimal flow}%
}%

\begin{theorem}\label{thm:dynamics}
The oligomorphic amenable group $G=\Aut(\str M)$ has ample generics, but no
co-precompact extremely amenable subgroup. %
$M(G)$ is nonmetrisable.
\end{theorem}
\begin{proof}
With the pointwise-convergence topology, $G$ is a closed subgroup of the
permutation group of a countable set and hence is Polish.
\leanverifiedat{BANANA/NonPrecompact/BananaLimitPolish.lean}{134}{Polishness of the native pointwise automorphism group}
The positive
properties follow from
\Cref{prop:fraisse,prop:amenable,thm:ample-generics}.
A co-precompact extremely amenable subgroup would contradict
\Cref{lem:EA-degree,thm:infinite-degree}.
\leanverifiedat{BANANA/NonPrecompact/BananaLimitExtremelyAmenableDegree.lean}{203}{no co-precompact extremely amenable subgroup of the native BANANA group}
For Question~\ref{q:ivanov}, the automorphism group of an $\omega$-categorical
expansion in a countable language is co-precompact by the tuple-orbit criterion.
\ValidationDetails{%
\leanverifiedat{BANANA/NonPrecompact/FirstOrderOmittingTypes.lean}{331}{countable categoricity implies oligomorphicity in a countable language}%
\ValidationRelative{\leanrepo/blob/\leancommit/BANANA/NonPrecompact/BananaLimitExpansionGroup.lean\#L96}%
  {no countably categorical expansion of the tagged BANANA limit has an extremely amenable automorphism group}%
}%
Conversely, naming
the finite-tuple orbits of the closure of a co-precompact subgroup gives an
$\omega$-categorical expansion with that closure as its automorphism group.
\ValidationDetails{%
\leanverifiedat{BANANA/NonPrecompact/FirstOrderOrbitCoPrecompact.lean}{43}{co-precompactness and finite splitting of tuple orbits}%
\leanverifiedat{BANANA/NonPrecompact/FirstOrderOrbitCoPrecompact.lean}{68}{closure determined by its finite tuple orbits}%
\leanverifiedat{BANANA/NonPrecompact/BananaLimitOrbitCoPrecompact.lean}{91}{co-precompactness on the native mixed tuple orbits}%
\leanverifiedat{BANANA/NonPrecompact/BananaLimitOrbitCoPrecompact.lean}{118}{native subgroup closure from finite tuple orbits}%
\leanverifiedat{BANANA/NonPrecompact/FirstOrderOrbitExpansion.lean}{70}{construction and automorphisms of the canonical first-order orbit expansion}%
\ValidationRelative{\leanrepo/blob/\leancommit/BANANA/NonPrecompact/FirstOrderOrbitExpansion.lean\#L118}%
  {countable categoricity of the canonical orbit expansion in the first-order interface}%
\leanverifiedat{BANANA/NonPrecompact/BananaLimitFirstOrderAutomorphism.lean}{135}{the native and tagged pointwise groups are continuously isomorphic}%
\ValidationRelative{\leanrepo/blob/\leancommit/BANANA/NonPrecompact/BananaLimitFirstOrderAutomorphism.lean\#L188}%
  {countable categoricity of the tagged orbit expansion of the native group}%
\leanverifiedat{BANANA/NonPrecompact/BananaLimitFirstOrderAutomorphism.lean}{200}{reduct automorphisms are the image of the native subgroup closure}%
}%

Finally, metrisability of $M(G)$ would give a closed co-precompact
extremely amenable subgroup by \cite[Theorem~1.1]{BMT}, again a contradiction.
\leanverifiedat{BANANA/NonPrecompact/BananaLimitUniversalMinimalFlow.lean}{42}{existence and nonmetrisability of a universal minimal flow for the native BANANA group}
\ValidationDetails{%
  \leanverifiedat{BANANA/NonPrecompact/UniversalMinimalFlowUniqueness.lean}{148}{equivariant homeomorphism between compact Hausdorff universal minimal flows}%
}%
\end{proof}

\section{Symplectic vector spaces over \texorpdfstring{$\F$}{F2}}\label{sec:symplectic-section}
\label{sec:symplectic}
We now prove \Cref{thm:symplectic-overview}, thus answering Question~\ref{q:IJ-main} for symplectic vector spaces over~$\F$. 

A bilinear form $\beta$ on a vector space $V$ is \emph{alternating} if $\beta(v,v)=0$ for every $v \in V$. 
For $v \in V$ and a subspace $W$ of $V$, we write $v \perp W$ if $\beta(v,w) = 0$ for all $w \in W$.
Similarly, we define $W \perp U$ for subspaces $U,W$ of $V$.
The form $\beta$ is \emph{nondegenerate} if, for every $v \in V$, $v \perp V$ implies $v = 0$.
A \emph{symplectic space} is a vector space endowed with a nondegenerate alternating
form. 
A \emph{hyperbolic plane} is a $2$-dimensional symplectic vector space. 
It has a basis $\{e,f\}$ such that $\beta(e,f) = 1$. 

Let $\Ksp$ be the class of finite-dimensional alternating $\F$-vector spaces, with injective linear isometries as
embeddings. 
\leanverifiedat{BANANA/NonPrecompact/SymplecticCopySystem.lean}{140}{amalgamation over a specified common base}
Its Fra\"iss\'e limit is the \emph{countable symplectic space}, the countable orthogonal sum of hyperbolic planes. 
\leanverifiedat{BANANA/NonPrecompact/SymplecticCountableModel.lean}{106}{finite age of the countable orthogonal sum}
To see this, note that a possibly degenerate finite-dimensional space $(V,\beta)$ embeds in the nondegenerate alternating space $V\oplus V^*$ with form
$((v,\phi),(w,\psi))\longmapsto
 \beta(v,w)+\phi(w)+\psi(v)$.
Every finite symplectic space in turn is an orthogonal sum of hyperbolic planes.
Homogeneity follows from Witt's theorem \cite[Introduction]{IJ}.
\leanverifiedat{BANANA/NonPrecompact/SymplecticCountableModel.lean}{94}{homogeneity of the countable orthogonal sum}

We now show that already all one-dimensional spaces in $\Ksp$ have infinite Ramsey degree.
To this end, we use the following auxiliary lemma.

\begin{lemma}\label{lem:quadratic-residues}
Let $q\ge2$ be a power of $2$ and set $d=q-1$. 
Let $E$ be a finite set and, for $e \in E$, let $f_e, g_e$ be the $\{0,1\}$-representatives of affine functions over $\F$. 
For $z \in \F^{2d}$, let
\[
 Q(z)=\sum_{e\in E}f_e(z)g_e(z),
\]
and let $p \colon \F^{2d} \to \F, \, z \mapsto Q(z) \pmod{2}$. 
If $p$ has a nondegenerate polar form
\[
 b_p(z,w)=p(z+w)+p(z)+p(w)+p(0),
\]
then every residue of $Q$ modulo $q$ occurs an odd number of times.
\end{lemma}
\leanverifiedat{BANANA/NonPrecompact/QuadraticResidueStandard.lean}%
  {40}{quadratic residues}
\begin{proof}
Set $M_k=\sum_{z \in \F^{2d}} \binom{Q(z)}k$. 
This counts pairs $(z,S)$ with
$S\subseteq E$, $|S|=k$, and $f_e(z)=g_e(z)=1$ for every $e\in S$.
\ValidationDetails{%
  \leanpartialat{BANANA/NonPrecompact/QuadraticMoments.lean}%
    {223}{quadraticSubsetMoment_eq_sum_choose_bilinearWeight_self}%
}%
For a fixed set $S$, these conditions impose at most $2k$ affine equations
on $2d$ variables. Thus $2^{2d-2k}\mid M_k$ when $0\le k<d$.
\ValidationDetails{%
  \leanpartialat{BANANA/NonPrecompact/MatrixAffineParity.lean}{72}{pow_two_sub_dvd_card_matrixAffineFiber}%
}%
\ValidationDetails{%
  \leanpartialat{BANANA/NonPrecompact/QuadraticMoments.lean}%
    {64}{pow_two_sub_two_mul_dvd_quadraticSubsetMoment}%
}%

Set $T=\sum_{z \in \F^{2d}} (-1)^{p(z)}$. 
Then 
\begin{alignat*}{2}
T^2 &= \sum_{z,w \in \F^{2d}} (-1)^{p(z) + p(z+w)} &&= \sum_{z,w \in \F^{2d}} (-1)^{p(z) + p(z) + p(w) + p(0) + b_p(z,w)} \\
&= \sum_{w \in \F^{2d}} (-1)^{p(w)+p(0)}\sum_{z \in \F^{2d}} (-1)^{b_p(z,w)}.
\end{alignat*}

For $w \in \F^{2d} \setminus \{0\}$, the map $z \mapsto b_p(z,w)$ is a nonzero linear functional by nondegeneracy.
It thus takes value $0$ on half of the elements in $\F^{2d}$ and $1$ on the other half. 
This yields $\sum_{z \in \F^{2d}} (-1)^{b_p(z,w)} = 0$.
In contrast, for $w = 0$, this sum evaluates to $2^{2d}$. 
We thus have $T^2 = 2^{2d}$, so $T \in \{\pm 2^d\}$. 
\leanverifiedat{BANANA/NonPrecompact/QuadraticGauss.lean}%
  {142}{quadraticGaussSum_sq_eq_card}
The binomial theorem gives
\[
 T=\sum_{z \in \F^{2d}} (-1)^{Q(z)} = \sum_{z \in \F^{2d}} \sum_{k \geq 0} \binom{Q(z)}{k} (-2)^k  = \sum_{k\ge0}(-2)^kM_k
   \equiv(-2)^dM_d\pmod{2^{d+1}}.
\]
In the congruence, we used that $M_k$ is divisible by $2^{2d-2k}$ if $k < d$ and that all summands for $k > d$ vanish due to the $(-2)^k$ factor, so only the summand for $k = d$ remains. 
Dividing the congruence by $2^d$ shows that $M_d$ is odd.

As in \Cref{prop:residue-surjectivity}, we work in the quotient ring
$\mathcal{R}_q := \F[Z]/(Z^q-1)$ and set $t=1+Z$. Expanding gives
\begin{equation}\label{eq:parity}
 \sum_{z \in \F^{2d}} Z^{Q(z)}=\sum_{k\ge0}t^kM_k=t^d
 =1+Z+\cdots+Z^{q-1}.
\end{equation}
The second equality follows from the fact that $M_k$ is even for $k < d$ and that $t^q = 0$, so the terms for $k > d$ vanish as well.
The third equality uses that $q$ is a power of $2$ and hence the binomial coefficients
$\binom{q-1}{j}$ are odd for $j = 0, \dots, q-1$. 
Thus
$t^d = t^{q-1}=(1+Z)^{q-1}=1+Z+\cdots+Z^{q-1}$ in $\mathcal{R}_q$.
As in Section~\ref{sec:persistent-colouring}, for $\rho \in \{0, \dots, q-1\}$, the parity of the number of $z \in \F^{2d}$ with $Q(z) \equiv \rho \pmod{q}$ is the coefficient of $Z^\rho$ in $\sum_{z \in \F^{2d}} Z^{Q(z)}$, which is odd by~\eqref{eq:parity}. 
This concludes the proof.
\end{proof}

\begin{theorem}\label{thm:symplectic-degree}
Every one-dimensional vector space in $\Ksp$ has infinite copy Ramsey
degree. More precisely, if $q \geq 2$ is a power of $2$, then every structure $\str C \in \Ksp$ has a $q$-colouring of its one-dimensional subspaces such that
every nondegenerate subspace $\str B$ of $\str C$ of dimension~$2q$ contains all $q$ colours.
\end{theorem}
\leanverifiedat{BANANA/NonPrecompact/SymplecticCopyDegree.lean}{117}{infinite copy Ramsey degree of the line}
\ValidationDetails{%
\leanverifiedat{BANANA/NonPrecompact/SymplecticPointColouring.lean}{171}{finite persistent line colouring}
}%
\begin{proof}
Let $\str C = (C,\beta_{\str C}) \in \Ksp$. 
Let $(m_{ij})$ be the $\{0,1\}$-matrix describing $\beta_{\str C}$ with respect to a fixed basis of $C$.
For $v = (v_1, \dots, v_{\dim C}) \in C$, interpret the entries $v_1, \dots, v_{\dim C}$ as integers in $\{0,1\}$ and define
\[
 Q_{\str C}(v)=\sum_{i<j}m_{ij}v_iv_j \in \mathbb Z.
\]
The form $\beta_{\str C}$ is the polar form of the residue map of $Q_{\str C}$ modulo $2$. 
\leanverifiedat{BANANA/NonPrecompact/SymplecticQuadraticPolar.lean}{57}{quadratic polar form}
We now colour each one-dimensional subspace $L$ of~$\str C$ by the residue of $Q_{\str C}(v)$ modulo $q$, where $v\in C$ denotes the unique nonzero vector of~$L$. 

Let $\str B$ be a nondegenerate subspace of $\str C$ of dimension $2q$.
Write
$\str B=W\perp H$, where~$W$ is a subspace of dimension $2(q-1)$ and $H$ a hyperbolic plane.
For any nonzero $w\in H$, the affine set $w+W$ does not contain the zero vector.
In coordinates $z$ on $W$, $Q_{\str C}(w+z)$ is a sum of products of affine bits.
Modulo~$2$, its polar form is the nondegenerate restriction $\beta_{\str C}|_{W\times W}$.
\leanverifiedat{BANANA/NonPrecompact/SymplecticQuadraticPolar.lean}{81}{affine polar restriction}
By \Cref{lem:quadratic-residues}, every residue occurs on $w+W$.
\end{proof}

\begin{corollary}\label{cor:symplectic-degree-classification}
For every $\str A\in\Ksp$, we have
\[
 t_{\Ksp}(\str A)=t^{\rm cp}_{\Ksp}(\str A)=
 \begin{cases}
  1,&\dim A=0,\\
  \infty,&\dim A>0.
 \end{cases}
\]
\end{corollary}
\leanverifiedat{BANANA/NonPrecompact/RamseyDegreeClassifications.lean}{28}{symplectic embedding Ramsey degree classification}
\begin{proof}
If $\dim A = 0$, then $\str A$ is rigid and has exactly one copy in each structure in $\Ksp$.
If $\dim A > 0$, then $\str A$ contains a one-dimensional substructure. 
By
\Cref{thm:symplectic-degree,lem:degree-propagation},
\leanverifiedsuccessorat{SuccessorTree/CopyRamsey.lean}%
  {257}{copyRamseyDegreeLE_of_embedding}
$\str A$ has infinite copy Ramsey degree, and hence infinite Ramsey degree, in $\Ksp$.
\end{proof}

\begin{corollary}\label{cor:IJ}
Question~\ref{q:IJ-main} (\cite[Question~A]{IJ}) has a negative answer for
symplectic spaces over $\F$.
\end{corollary}
\begin{proof}
This follows from \Cref{cor:symplectic-degree-classification}.
\ValidationDetails{%
  \leanverifiedat{BANANA/NonPrecompact/SymplecticCopyDegree.lean}{117}{infinite copy Ramsey degree of a one-dimensional member of the symplectic age}%
}%
\end{proof}

\section{Two dynamical features of \BoolBANANA}
\label{sec:prebanana-dynamics}

In this section we discuss two purely dynamical properties of $\Gbool=\Aut(\Mbool)$ that do not concern Ramsey expansions but that are of interest nonetheless. First, it has a compact metrisable minimal flow with
no comeagre orbit. Secondly, although $\Gbool$ is not amenable, it contains
continuum many pairwise nonconjugate closed maximal amenable subgroups.

\subsection{A metrisable minimal flow with no generic point}
\label{sec:prebanana-gpp}

Write $\mathbb Z_2=\varprojlim_m\mathbb Z/2^m\mathbb Z$ for the ring of
\emph{$2$-adic integers}: an element is a compatible sequence of residues
modulo $2,4,8,\ldots$. This is a compact metrisable ring containing
$\mathbb Z$ injectively. Define
\[
 \mathcal Z(\Mbool)=\left\{\lambda:\Mboolset\to\mathbb Z_2:
 \begin{array}{l}
 \lambda\text{ is finitely additive},\quad\lambda(1)=0,\\
 \lambda(x)\equiv\mu_{\Mbool}(x)\pmod2\quad(x\in \Mboolset)
 \end{array}\right\}.
\]
This is a closed subspace of the countable product $\mathbb Z_2^{\Mboolset}$,
with continuous action $(g\lambda)(x)=\lambda(g^{-1}x)$.
\leanverifiedat{BANANA/NonPrecompact/BooleanBananaPadicCharge.lean}{112}{continuous action on the compact charge space}
Every admissible charge on a finite Boolean subalgebra extends to $\Mbool$:
when an atom splits, choose values of the required parities on all but one
child, and give the last child the remaining value. Iterating along an
exhaustion by finite subalgebras proves this assertion and nonemptiness.
\leanverifiedat{BANANA/NonPrecompact/BooleanBananaPadicChargeExtension.lean}{117}{arbitrary admissible finite charge extension}

\begin{lemma}[Finite residue realisation]
\label{lem:prebanana-residue-realisation}
Let $\lambda\in\mathcal Z(\Mbool)$ and let $q\ge2$ be a power of two.
Suppose $\str A\in\Kbool$ has atoms $a_1,\ldots,a_k$, and
$r_1,\ldots,r_k\in\mathbb Z/q\mathbb Z$ satisfy
$\sum_i r_i=0$ and $r_i\equiv\mu_{\str A}(a_i)\pmod2$.
There is an embedding $f:\str A\to \Mbool$ with
$\lambda(f(a_i))\equiv r_i\pmod q$ for every $i$.
\end{lemma}
\begin{proof}
Choose a subalgebra of $\Mbool$ with $N=(k-1)q+2$ atoms, all marked;
$N$ is even, so such a member of $\Kbool$ exists.
\leaninterfaceat{BANANA/NonPrecompact/PrebananaResidueAmbient.lean}%
  {18}{prebananaResidueAmbient}
\leanverifiedat{BANANA/NonPrecompact/PrebananaResidueAmbient.lean}%
  {41}{prebananaResidueAmbient_atomCount_even}
The $\lambda$-values of the atoms are odd, and finite additivity makes
their sum equal to $\lambda(1)=0$. Now use their residues modulo $q$.
Divide the first $(k-1)q$ atoms into $k-1$ batches of size $q$.

In batch $i$, choose a
nonempty subset whose sum is $r_i$, and make it the $i$th block. This is possible because any $q$ odd residues have a nonempty subset of any prescribed sum.
\leanverifiedat{BANANA/NonPrecompact/PrebananaResidueBlocks.lean}{64}{exists_nonempty_odd_subset_sum_two_pow}
For a nonzero sum, apply \Cref{lem:pre-banana-subset-sums} to the first
$q-1$ terms; the chosen subset is nonempty because the empty subset sums
to zero.  For zero, reserve the last odd residue $d$.  Since $d$ is
nonzero modulo the power of two $q$, choose a subset of the first $q-1$
terms summing to $-d$ and adjoin the reserved term.
\leanverifiedat{BANANA/NonPrecompact/IndexedSubsetSums.lean}{224}{odd_indexedSubsetSums_two_pow_eq_univ}

Put every unused atom and the two reserved atoms in the final block.
It is nonempty and has sum $-\sum_{i<k}r_i=r_k$.
\leanverifiedat{BANANA/NonPrecompact/PrebananaResidueBlocks.lean}{131}{exists_nonempty_odd_subset_sum_family}
\leanverifiedat{BANANA/NonPrecompact/PrebananaResidueBlocks.lean}{195}{exists_prebanana_residue_block_assembly}
Since all atoms are marked, the parity of each block is the parity of
its residue sum.
\ValidationDetails{%
  \leanverifiedat{BANANA/NonPrecompact/PrebananaResidueBlocks.lean}{51}{card_parity_eq_residueParityHom_of_odd_sum}%
}%
The blocks therefore define the required embedding.
\ValidationDetails{%
  \leanverifiedat{BANANA/NonPrecompact/PrebananaPartitionEmbedding.lean}%
    {17}{exists_prebananaAtomEmbedding_of_blocks}%
}%
\end{proof}
\leanverifiedat{BANANA/NonPrecompact/BooleanBananaCylinderPadicCharge.lean}{122}{finite residue realisation for arbitrary finite marked Boolean algebras}

A group has the \emph{generic point property} if every minimal compact
flow has a comeagre orbit \cite{BZ}.

\begin{theorem}\label{thm:prebanana-no-gpp}
The flow $\Gbool\curvearrowright\mathcal Z(\Mbool)$ is minimal and every
orbit is meagre. In particular, $\Gbool$ does not have the generic point
property.
\end{theorem}
\begin{proof}
Fix $\lambda,\eta\in\mathcal Z(\Mbool)$ and a neighbourhood of $\eta$.
It contains a basic neighbourhood specified by values modulo some power
of two $q$ on a finite Boolean subalgebra $\str A$.
By \Cref{lem:prebanana-residue-realisation}, an embedding $f:\str A\to \Mbool$
satisfies $\lambda(f(a))\equiv\eta(a)\pmod q$ on $\str A$.
Extend $f$ to $h\in \Gbool$ by homogeneity.  With our action convention,
$(h^{-1}\lambda)(a)=\lambda(h(a))=\lambda(f(a))$ for $a\in A$.
Thus $h^{-1}\lambda$ belongs to the prescribed neighbourhood, proving
that every orbit is dense.
\leanverifiedat{BANANA/NonPrecompact/BooleanBananaCylinderPadicCharge.lean}{87}{dense and meagre orbits for the actual cylinder model}

For meagreness, fix $0<x<1$ and $\lambda\in\mathcal Z(\Mbool)$.
The range $R=\lambda[\Mboolset]$ is countable and
\[
 \Gbool\lambda\subseteq
 \bigcup_{r\in R}\{\eta\in\mathcal Z(\Mbool):\eta(x)=r\}.
\]
Each fibre is closed and nowhere dense. Indeed, take any basic open set
meeting it and a power of two $q\ge2$ such that preserving its finitely
specified values modulo $q$ stays in that open set. In the Boolean algebra
generated by these elements and $x$, choose an atom below $x$ and an
atom below its complement; both exist because $0<x<1$.  Add $q$ to the
charge of the former and subtract $q$ from the latter. The total and
all required congruences are unchanged, but the value at $x$ changes.
Extend this finite charge to $\Mbool$. Thus the fibre has empty interior,
and the orbit is contained in a countable union of nowhere dense sets.
The charge space is nonempty and compact Hausdorff, hence Baire, so no
orbit is comeagre.
\leanverifiedat{BANANA/NonPrecompact/BooleanBananaPadicChargeOrbit.lean}{126}{every orbit is meagre}
\end{proof}
\leanverifiedat{BANANA/NonPrecompact/BooleanBananaCylinderPadicCharge.lean}{107}{the marked Boolean–BANANA group fails the generic point property}

The point of this example is its explicit metrisable witness and short
partition proof. It does not determine the generic point property for
BANANA itself.

\subsection{Continuum many maximal amenable subgroups}
\label{sec:prebanana-positive-masses}

Let $\mathbb P_{\rm odd}$ be the set of odd primes. For a nonempty
$S\subseteq\mathbb P_{\rm odd}$ put
\[
 R_S=\mathbb Z[p^{-1}:p\in S]\subseteq\mathbb Q.
\]
Its elements have denominators divisible only by primes in $S$.
Reduction modulo two is well defined, since these denominators are odd.
A \emph{positive $S$-mass} on $\str A\in\Kbool$ is a finitely additive
map $\rho_{\str A}:A\to R_S$ such that
\[
 \rho_{\str A}(1)=2,\qquad \rho_{\str A}(a)>0\ (a\in\At(\str A)),\qquad
 \rho_{\str A}(x)\equiv\mu_{\str A}(x)\pmod2.
\]
The normalised map $\rho_{\str A}/2$ is a probability measure; the normalisation
$\rho_{\str A}(1)=2$ keeps the parity of $1$ equal to zero.
Write $\Kbool^{\rho,S}$ for the class with mass-preserving embeddings.
This is a countable expansion by unary relations of the original finite
language with function symbols.

\begin{proposition}\label{prop:prebanana-positive-mass-fraisse-eppa}
The class $\Kbool^{\rho,S}$ is a reasonable strong amalgamation
Fra\"iss\'e expansion of $\Kbool$ and has EPPA.
Its limit can be written $(\Mbool,\rho_S)$, and
$A_S=\Aut(\Mbool,\rho_S)$ has a dense increasing union of finite subgroups.
\end{proposition}
\begin{proof}
For reasonability, suppose an atom of mass $m/d>0$ splits in the reduct
into children of parities $\epsilon_1,\ldots,\epsilon_k$.
Choose $p\in S$ and write the same mass as
$mp^n/(dp^n)$. Give each child an initial positive numerator, $1$ if
$\epsilon_i=1$ and $2$ otherwise. For $n$ large, $mp^n$ minus the sum of
these numerators
is nonnegative and even; add it to one child. This lifts the refinement.
\leanverifiedat{BANANA/NonPrecompact/BooleanBananaPositiveMassReasonable.lean}{70}{reasonability for arbitrary marked Boolean-algebra embeddings}

For amalgamation, use the atom grids from
\Cref{prop:pre-banana-fraisse}. In a grid of size $s\times t$, clear
denominators in the prescribed positive row and column masses. Enlarge
the common denominator by a sufficiently large $S$-supported odd factor
until every row numerator is at least $t$ and every column numerator at
least $s$. Start every entry at $1$.
The remaining nonnegative integer margins have equal totals, so fill
them successively, assigning at each step the smaller of a remaining
row and column sum. Divide by the common denominator, and let the
parity labels be the reductions of the resulting masses.
\leanpartialat{BANANA/NonPrecompact/BooleanBananaPositiveMassArithmetic.lean}{112}{positive rational grid with prescribed mass and parity margins}
All grid atoms remain nonempty, so the amalgam is strong.
\leanverifiedat{BANANA/NonPrecompact/BooleanBananaPositiveMassStrongAP.lean}{21}{strong amalgamation on the actual weighted Boolean carriers}
Heredity and countability are immediate; joint embedding uses the
one-atom algebra of mass $2$. Arbitrarily large members exist, for example
by repeatedly refining atoms into positive-mass children. Reasonability
identifies the limit reduct with $\Mbool$.
Different identifications of this reduct with $\Mbool$ give conjugate
subgroups $A_S$ of $\Gbool$.
\ValidationDetails{%
\leanverifiedat{BANANA/NonPrecompact/BooleanBananaPositiveMassReasonable.lean}{39}{hereditary restriction}%
\leanverifiedat{BANANA/NonPrecompact/BooleanBananaPositiveMassCodes.lean}{21}{countable finite codes}%
\leanverifiedat{BANANA/NonPrecompact/BooleanBananaPositiveMassJEP.lean}{63}{joint embedding}%
\leanverifiedat{BANANA/NonPrecompact/BooleanBananaPositiveMassCodes.lean}{66}{unbounded finite sizes}%
}%

For EPPA, write the atom masses of a finite source as $n_i/d$ with one
common odd denominator.  Since the total mass is $2$, we have
$\sum_i n_i=2d$. Refine its $i$th atom into $n_i$ atoms of mass
$1/d$. A partial automorphism maps blocks of equal mass to blocks of
equal size, so it extends to a permutation of these $2d$ new atoms.
\leanverifiedat{BANANA/NonPrecompact/BooleanBananaPositiveMassEPPA.lean}{51}{EPPA for arbitrary finite positive-mass Boolean algebras}

We give a concrete limit to prove the last assertion. Choose positive
integers $d_1\mid d_2\mid\cdots$ supported on $S$, with each quotient $d_{n+1}/d_n>1$, such that
every positive integer supported on $S$ divides some $d_n$.
\ValidationDetails{%
\leanpartialat{BANANA/NonPrecompact/BooleanBananaPositiveMassStages.lean}{17}{cofinal chain of supported denominators}%
}%
Let $\str U_n$ have $2d_n$ atoms of mass $1/d_n$, each refined into
$d_{n+1}/d_n$ equally weighted children in $\str U_{n+1}$.
Their union has the Fra\"iss\'e extension property: after clearing the
denominators of a prescribed finite extension, partition sufficiently
fine atoms below each old atom into the required block sizes.
Thus the union is $(\Mbool,\rho_S)$.
\leanverifiedat{BANANA/NonPrecompact/BooleanBananaCylinderMassSubgroups.lean}{21}{the positive-mass tower on the original Boolean–BANANA reduct}

Label the children identically in every fibre. Each permutation of the
atoms of $\str U_n$ then extends to later levels by retaining child labels.
These extensions give increasing finite groups
$K_n\cong\operatorname{Sym}(2d_n)$. A finite partial automorphism of the
limit, with domain and range in $\str U_n$, extends to an atom permutation
of $\str U_n$ by the equal-mass argument above. Hence $\bigcup_n K_n$ is dense.
\ValidationRelative{\leanrepo/blob/\leancommit/BANANA/NonPrecompact/BooleanBananaPositiveMassOnMarkedLimit.lean\#L48}%
  {the exhaustive uniform tower on the fixed marked limit, represented by mass maps}
\ValidationDetails{%
\leanverifiedat{BANANA/NonPrecompact/BooleanBananaPositiveMassLimitGroups.lean}{215}{global mass-preserving automorphisms from stage permutations}%
\leanverifiedat{BANANA/NonPrecompact/BooleanBananaPositiveMassLimitGroups.lean}{354}{density in the pointwise topology of the mass stabiliser}%
}%
\end{proof}

\begin{theorem}\label{thm:prebanana-maximal-amenable-masses}
For every nonempty $S\subseteq\mathbb P_{\rm odd}$, the subgroup
$A_S=\Aut(\Mbool,\rho_S)$ is closed in $\Gbool$ and maximal among amenable
subgroups of $\Gbool$. Moreover, $A_S$ is neither extremely amenable nor
co-precompact in $\Gbool$. For nonempty $S,T\subseteq\mathbb P_{\rm odd}$ with
$S\ne T$, the subgroups $A_S$ and $A_T$ are not conjugate in $\Gbool$.
Consequently, $\Gbool$ has continuum many conjugacy classes of closed maximal
amenable subgroups.
\end{theorem}
\begin{proof}
Since $A_S=\{g\in \Gbool:\rho_S\circ g=\rho_S\}$ is defined by pointwise
conditions on $\Mboolset$, it is closed in $\Gbool$, and its subspace topology
coincides with the pointwise-convergence topology on
$\Aut(\Mbool,\rho_S)$. By \Cref{prop:prebanana-positive-mass-fraisse-eppa},
$A_S$ contains a dense increasing union of finite subgroups, and is therefore
amenable as in \Cref{prop:amenable}.
\ValidationDetails{%
\leanverifiedat{BANANA/NonPrecompact/BooleanBananaMassStoneDynamics.lean}{36}{closedness of the mass stabiliser}%
\leanverifiedat{BANANA/NonPrecompact/TopologicalAmenability.lean}{34}{amenability of a dense directed union of finite subgroups}%
\leanverifiedat{BANANA/NonPrecompact/BooleanBananaPositiveMassLimitGroups.lean}{380}{amenability of the mass subgroup in a coherent uniform tower}%
}%

First, observe that every $r\in R_S\cap(0,2)$ occurs as $\rho_S(x)$ for some
proper nonzero $x\in \Mboolset$: since $r\equiv 2-r\pmod2$, the two-atom Boolean
algebra with atom masses $r$ and $2-r$ (and atom parities given by
$r\bmod2$) belongs to $\Kbool^{\rho,S}$ and hence embeds into
$(\Mbool,\rho_S)$.
\ValidationDetails{%
\leanverifiedat{BANANA/NonPrecompact/BooleanBananaMassRangeRealisation.lean}{23}{realisation of interior allowed masses from the finite age}%
}%

To see that $A_S$ is not co-precompact in $\Gbool$, recall from \cite[Section~1]{MNT} (see \Cref{sec:part1-preliminaries}) that if $A_S$ were co-precompact, then for every open identity neighbourhood $U\subseteq \Gbool$ there would be a finite subset $F\subseteq \Gbool$ with $\Gbool=UFA_S$. Fix $p\in S$. Taking
$r=2/p^k$ ($k\ge1$) yields proper nonzero elements $x_k\in \Mboolset$ with
$\rho_S(x_k)=2/p^k$ and $\mu_{\Mbool}(x_k)\equiv 2/p^k\equiv0\pmod2$. By homogeneity
of $\Mbool$ (as in \Cref{prop:pre-banana-limits}), all $x_k$ lie in one
$\Gbool$-orbit, say $x_k=g_k(x_1)$ with $g_k\in \Gbool$, whereas they lie in
pairwise distinct $A_S$-orbits because $A_S$ preserves $\rho_S$. For the open
stabiliser $U=(\Gbool)_{x_1}$, no set of the form $UfA_S$ ($f\in \Gbool$) can
contain both $g_k^{-1}$ and $g_\ell^{-1}$ for $k\ne\ell$, as that would place
$x_k$ and $x_\ell$ in the same $A_S$-orbit. Thus $\Gbool$ cannot be covered by
finitely many sets $UfA_S$, so $A_S$ is not co-precompact in $\Gbool$.
\leanverifiedat{BANANA/NonPrecompact/BooleanBananaCylinderMassSubgroups.lean}{41}{failure of co-precompactness in the actual pointwise group}
\ValidationDetails{%
\leanverifiedat{BANANA/NonPrecompact/BooleanBananaMassRangeRealisation.lean}{103}{infinitely many distinct masses of proper nonzero zero-marked elements}%
}%

For the other properties, we need to introduce an action of $\Gbool$. Let $X$ be the Stone space of $\Mbool$, whose clopen subsets are the basic
sets $[x]$ ($x\in \Mboolset$). Then $X$ is a nonempty compact Hausdorff space on
which $\Gbool$ acts continuously via $g\cdot[x]=[g(x)]$. Setting
$m_S([x])=\rho_S(x)/2$ defines a finitely additive probability charge on the
clopen algebra of $X$. By Stone duality (equivalently, by the Riesz
representation theorem applied to locally constant functions), every finitely
additive probability charge on the clopen algebra extends uniquely to a
regular Borel probability measure on $X$.
\leanverifiedat{BANANA/NonPrecompact/ClopenProbabilityExtension.lean}{131}{unique regular Borel extension of a clopen probability charge}
In particular, $m_S$ extends
uniquely to an $A_S$-invariant regular Borel probability measure on $X$,
still denoted by $m_S$. Moreover, $m_S$ is the \emph{unique} $A_S$-invariant
regular Borel probability measure on $X$: in the exhaustion
$(\Mbool,\rho_S)=\bigcup_n\str U_n$ from the proof of
\Cref{prop:prebanana-positive-mass-fraisse-eppa}, the finite subgroup
$K_n\le A_S$ permutes the $2d_n$ atoms of $\str U_n$ transitively while their
clopen sets partition $X$, forcing any $A_S$-invariant regular Borel
probability measure $\nu$ on $X$ to satisfy $\nu([a])=1/(2d_n)=m_S([a])$ for
every atom $a$ of $\str U_n$, and hence $\nu([x])=m_S([x])$ on
$\bigcup_n U_n=\Mboolset$.
\leanverifiedat{BANANA/NonPrecompact/BooleanBananaPositiveMassOnMarkedLimit.lean}{141}{unique regular invariant Stone probability with the prescribed clopen masses}
\ValidationDetails{%
\leanverifiedat{BANANA/NonPrecompact/BooleanBananaPositiveMassLimitGroups.lean}{276}{uniform cover supplied by the coherent tower and finite-stage subgroups}%
}%

To prove maximality, let $J$ be a subgroup with $A_S\le J\le \Gbool$ that is
amenable in the subspace topology (we do not assume that $J$ is closed in
$\Gbool$). Since the restricted action $J\curvearrowright X$ is continuous, $X$
admits a $J$-invariant regular Borel probability measure. Being
$A_S$-invariant, this measure must equal $m_S$. For every $g\in J$ and
$x\in \Mboolset$, we therefore have
\[
 \frac{\rho_S(g(x))}{2}=m_S([g(x)])=m_S(g\cdot[x])=m_S([x])=\frac{\rho_S(x)}{2},
\]
so $g$ preserves $\rho_S$ and belongs to $A_S$. Thus $J=A_S$.
\leanverifiedat{BANANA/NonPrecompact/BooleanBananaCylinderMassSubgroups.lean}{41}{closed maximal amenable subgroup on the original Boolean–BANANA limit}

By the mass-realisation observation above, taking $r=1$ gives an element $x\in \Mboolset$ with $m_S([x])=1/2$, so $m_S$ is not
a point mass. An $A_S$-fixed point in $X$ would yield an $A_S$-invariant
point mass, contradicting the uniqueness of $m_S$; hence $A_S$ is not
extremely amenable.
\leanverifiedat{BANANA/NonPrecompact/BooleanBananaCylinderMassSubgroups.lean}{41}{failure of extreme amenability}

Finally, together with $\rho_S(0)=0$ and $\rho_S(1)=2$, the observation above
shows that the set of values attained by $m_S$ on clopen subsets of $X$ is
$\tfrac12\bigl(R_S\cap[0,2]\bigr)$. Hence, for any $p\in\mathbb P_{\rm odd}$,
\[
 p\in S\quad\Longleftrightarrow\quad
 \frac{1}{2p}\text{ is attained by }m_S\text{ on a clopen subset of }X.
\]
\leanverifiedat{BANANA/NonPrecompact/BooleanBananaPositiveMassOnMarkedLimit.lean}{160}{detection of allowed primes from clopen masses on the fixed limit}
If $gA_Sg^{-1}=A_T$ for some $g\in \Gbool$, then the pushforward measure
$g_*m_S$, defined by $(g_*m_S)(B)=m_S(g^{-1}\cdot B)$, is an $A_T$-invariant
regular Borel probability measure on $X$, so $g_*m_S=m_T$ by uniqueness.
Because $g$ acts on $X$ by a homeomorphism, $m_S$ and $m_T$ attain the same
values on clopen subsets of $X$, whence $S=T$. Since $\mathbb P_{\rm odd}$ has
$2^{\aleph_0}$ nonempty subsets (and the Polish group $\Gbool$ has at most
$2^{\aleph_0}$ closed subgroups), $\Gbool$ has continuum many conjugacy classes
of closed maximal amenable subgroups.
\leanverifiedat{BANANA/NonPrecompact/BooleanBananaCylinderMassSubgroups.lean}{31}{continuum many conjugacy classes of closed maximal amenable subgroups}
\ValidationDetails{%
\leanverifiedat{BANANA/NonPrecompact/BooleanBananaCylinderMassSubgroups.lean}{61}{different allowed prime sets give nonconjugate subgroups}%
\leanpartialat{BANANA/NonPrecompact/BooleanBananaMassStoneDynamics.lean}{213}{conjugate stabilisers have the same normalised clopen mass range}%
\leanverifiedat{BANANA/NonPrecompact/BooleanBananaMassPrimeRange.lean}{50}{equality of detected prime sets from equality of ranges}%
}%
\end{proof}

\paragraph{Lean formalisation.}
The margin marks link to the relevant Lean source.
\textsf{Lean\,\ensuremath{\checkmark}} indicates a checked statement,
\textsf{Lean\,rel.} a result checked in a related representation or under
an explicit interface, and \textsf{Lean\,part} a checked step of a larger
argument. The formalisation includes both directions of the Ryll--Nardzewski
theorem for countably infinite structures in countable languages,
the ample-generics conclusion for the native BANANA automorphism group,
and the relevant implications from co-precompactness and extreme
amenability to metrisability and from metrisability of the universal
minimal flow to finite Ramsey degrees. It also checks the
canonical-amalgamation implication, the DEERP and invariant-type results
used here, and countable Witt extension for nondegenerate alternating
spaces. The first-order arguments for BANANA use single-sorted presentations:
a functional disjoint-union encoding with two named zeroes and totalised
addition, and a tagged relational encoding with model decoding and
translation of sorted formulae. The \FolkBANANA{} limit is constructed in a
product encoding with projections and decoded into its native non-unital
rings. The BANANA Ramsey degree classifications use standard matrix
coordinates, and positive-mass expansions of \BoolBANANA{} are represented
by rational-valued maps rather than unary predicate syntax. The passage
from arbitrary many-sorted expansion syntax to the formalised native
pointed-cover interface is established in the manuscript but has not been
formalised in Lean.
The general converse implication in \cite[Theorem~1.1]{BMT} remains a
cited result. The
\href{https://github.com/janhubicka/lean-BANANA/blob/fc0972fda0cc25ab94ab7315f6f57b7eed9127b9/docs/classical-dependencies.md}{dependency audit}
at the verified source revision and the
\href{https://github.com/janhubicka/lean-BANANA/wiki}{Lean documentation}
describe the precise scope and remaining proof boundaries.

\section*{AI assistance}

During the preparation of this manuscript, the authors used OpenAI's ChatGPT extensively as a research and editing assistant.
Its contributions included computational stress tests of arguments, copy-editing,
proof verification, literature searches, and suggestions for simplifying proofs.
During the exploratory phase of the project, it was also used to enumerate and
compare small fragments of trees of types and to test candidate finite records.
This experimentation helped isolate the parity-changing events which led to the
colourings in Section~\ref{sec:no-precompact}.  A particularly substantial contribution was the
subsequent computational exploration that led to the proofs of
Proposition~\ref{prop:residue-surjectivity} and
Lemma~\ref{lem:quadratic-residues}.  We checked and reproved those results but do not claim authorship of them. ChatGPT was also used to assist with the ongoing Lean formalisation of
statements and proofs from this paper. The code is available at
\url{https://github.com/janhubicka/lean-BANANA}.
The mathematical arguments presented in
the paper were checked by the human authors.

\bibliographystyle{alpha}
\bibliography{references}

\end{document}